\documentclass[]{aiaa-tc}% insert '[draft]' option to show overfull boxes
\pdfoutput=1

 \title{Generalised Lebesgue Stable Flux Reconstruction}

 \author{
  Will Trojak%
    \thanks{PhD Candidate, wt247@cam.ac.uk, Department of Engineering, University of Cambridge, Cambridge, UK, CB2 1PZ}
    }

 \AIAApapernumber{YEAR-NUMBER}
 \AIAAconference{Conference Name, Date, and Location}
 \AIAAcopyright{\AIAAcopyrightD{YEAR}}

 \newcommand{\etal}{\emph{et~al.}}
 \newcommand\floor[1]{\lfloor#1\rfloor}

 \newcommand{\hfd}{\hat{f}^{\delta}}
 
 \newcommand{\hfI}{\hat{f}^{I}}
 \newcommand{\hfC}{\hat{f}^{C}}
 \newcommand{\hud}{\hat{u}^{\delta}}
 \newcommand{\hu}{\hat{u}}
 \newcommand{\huI}{\hat{u}^{I}}

 \newcommand\px[2]{\frac{\partial #1}{\partial {#2}}}
 \newcommand\pxi[3]{\frac{\partial^{#1}#2}{\partial {#3}^{#1}}}
 \newcommand\dx[2]{\frac{d #1}{d #2}}
 \newcommand\dxi[3]{\frac{d^{#1}#2}{d {#3}^{#1}}}

\graphicspath{ {Figs/} }

\usepackage{graphicx}
\usepackage{graphics}
\usepackage{upgreek}
\usepackage{float}
\usepackage{epstopdf}
\usepackage{amsmath}
\usepackage{wrapfig}
\usepackage{amssymb}
\usepackage{pdflscape}
\usepackage{listings}
\usepackage{cite}
\usepackage{subcaption}
\usepackage{multirow}
\usepackage[percent]{overpic}
\usepackage{multicol}
\usepackage{multirow}
\usepackage{mathtools}
\usepackage{color}
\usepackage{setspace}
\usepackage{algorithm}
\usepackage{algorithmic}

\begin{document}

\maketitle

\begin{abstract}
	A unique set of correction functions for Flux Reconstruction is presented, with there derivation stemming from proving the existence of energy stability in the Lebesgue norm. The set is shown to be incredibly arbitrary with the only union to existing correction function sets being show to be for Nodal DG via reconstruction. Von Neumann analysis of both linear advection and coupled advection-diffusion is used to show that once coupled to a temporal integration method, good CFL performance can be achieved and correction functions presented that may have better dispersion and dissipation for application to implicit LES. Lastly, the turbulent Taylor-Green vortex test case is then used to show that correction functions can be found that improve the accuracy of the scheme when compared to the error levels of Discontinuous Galerkin.
\end{abstract}

\section{Introduction}
	Flux reconstruction as a doctrine for solving hyperbolic systems of equations has existed in various forms for considerable time, for example the now ubiquitous ENO~\cite{Harten1987} and WENO~\cite{Liu1994} methods. However, in recent years more advanced methods of flux reconstruction have come to the for, initially through Discontinuous Galerkin (DG)~\cite{Reed1973,Cockburn2000} and then Spectral Element (SE)~\cite{Karniadakis2013} and Spectral Volume (SV)~\cite{Wang2002} methods. These methods introduce a concept that would become key in the development of what we now recognise as the Flux Reconstruction (FR). Namely the use of domain sub-division to allow local polynomial fitting, with the motivation being simpler handling of unstructured grids and high-order. FR can now be defined as a method that utilises sub-domains and polynomial correction to formulate a piecewise continuous solution as was introduced by Correction Procedure via Reconstruction (CPR)~\cite{Gao2013a,Wang2009} together with the works of Huynh~\cite{Huynh2007,Huynh2009}. Extensions have since been made to generalise the topology of solutions to simplicies and hypercubes~\cite{Williams2011,Huynh2011,Castonguay2012a,Williams2014,Sheshadri2016a} to reasonable effect, although some questions still remain.
	
	The focus of this research is, however, not concerned with topology or even with a particular equation set, but with the correction procedure definition. This is of particular concern, as was shown by Vincent~\etal~\cite{Vincent2011} and Vermeire~and~Vincent~\cite{Vermeire2016}, the correction function has enormous control over the character of the numerical scheme. Yet, there are still some established correction function, for example the Lumped-Chebyshev-Lobatto correction of Huynh~\cite{Huynh2007}, which have not yet been shown to fit within the current realm of stable schemes defined by there correction. Hence, our current view of FR schemes must be incomplete. The current families of corrections will be introduced in the mathematical description of FR, however for reference they may be found in \cite{Vincent2010,Vincent2015,Trojak2018}.
	
	The investigation presented in this paper will have the following structure. Firstly, the method of FR will be presented together with the status-quo of FR correction functions. Secondly, we will present an energy stability proof that simplifies the Sobolev stability proof into the Lebesgue space, which is then shown to give rise to a new and expansive set of corrections. Thirdly, we show that this arbitrary form of FR does not contradict mass conservation of the underlying FR method and lastly we will provide analysis, through both advective-diffusive von Neumann analysis and numerical tests on the turbulent incompressible Taylor-Green vortex.

\section{Flux Reconstruction} 
	To demonstrate the mathematical procedure of Flux Reconstruction~\cite{Huynh2007,Vincent2010,Jameson2012} we will consider a 1D conservation equation on the domain $\mathbf{\Omega}$ which is divided into sub-domains such that:
	\begin{equation}
		\mathbf{\Omega} = \bigcup^N_{i=0} \mathbf{\Omega}_i \quad \mathrm{and} \quad \mathbf{\Omega}_i \cap \mathbf{\Omega}_j = \emptyset \:\: \forall \: i\ne j
	\end{equation}
	we then define the reference domain for which the sub-domains may be project into as $\hat{\mathbf{\Omega}}$. With the projection defined such that:
	\begin{equation}
		\hat{u}_j = J_ju_j
	\end{equation}
	where: $u$ is a conserved variable; $J_j$ is the $j^{\mathrm{th}}$ element Jacobian; and hats and used to signifies a variable in the reference domain. With the domain and decomposition set out, we will now present 1D conservation equation to be considered:
	\begin{equation}\label{eq:conservation}
		\px{u}{t} + \px{f}{x} = 0
	\end{equation}
	where $f$ is some flux function. Within each sub-domain a discontinuous $p^{\mathrm{th}}$ order polynomial can be constructed from data stored at the $p+1$ solution points within the sub-domain, such that:
	\begin{align}
		\hud &= \sum_{j=0}^p \hud(\xi_j)l_j(\xi) \\
		\hfd &= \sum_{j=0}^p \hfd(\xi_j)l_j(\xi)
	\end{align}
	where $\delta$ is used to indicate that a variable is discontinuous and $l_j(\xi)$ is the $j^{\mathrm{th}}$ Lagrange basis polynomial for $\xi \in \hat{\mathbf{\Omega}}$, defined by:
	\begin{equation}
		l_j(\xi) = \prod_{i=0}^p \bigg(\frac{\xi-\xi_i}{\xi_j-\xi_i}\bigg)(1-\delta_{ij})
	\end{equation}
	
	Approximation of the solution to Eq.(\ref{eq:conservation}) requires that the domain wide solution is $C^0$ continuous and to fulfil this requirement a correction procedure is followed. The first step of which is to use the discontinuous flux polynomial to interpolate the flux function to collocated flux points on the boundary of each sub-domain. In the 1D case this gives us $\hfd_j(-1)=\hfd_{j,L}$ and $\hfd_j(1)=\hfd_{j,R}$. Once, the interpolated values at the collocated edges points, or flux points, are calculated, an appropriate Riemann solver~\cite{Toro2009} may used to calculate a common interface flux value, $\hfI_{j,L}$ and $\hfI_{j,R}$ for the left and right interfaces respectively. This then allows a correction to be performed by propagating the correction to the shared common interface value in the adjacent elements to form the continuous flux polynomial $\hfC$. This is done in the following manner:
	\begin{equation}
		\hfC_j(\xi) = (\hfI_{j,L}-\hfd_{j,L})h_L(\xi) + (\hfI_{j,R} - \hfd_{j,R})h_R(\xi) + \hfd_j
	\end{equation}
	where $h_L$ and $h_R$ are the left and right correction function respectively. With the boundary conditions:
	\begin{align}
		h_L(-1) = 1 \quad &\mathrm{and} \quad h_L(1) = 0 \\
		h_R(-1) = 0 \quad &\mathrm{and} \quad h_R(1) = 1
	\end{align}
	Upon differentiation of $\hfC$ the 1D conservative equation can be defined as:
	\begin{equation}
		\px{\hu_j}{t} = -\sum^p_{i=0}\hfd_{j}(\xi_i)\dx{l_i(\xi)}{\xi} - (\hfI_{j,L}-\hfd_{j,L})\dx{h_L(\xi)}{\xi} - (\hfI_{j,R}-\hfd_{j,R})\dx{h_R}{\xi} \label{eq:FR_broken}
	\end{equation}
	Time marching of the solution can be performed via an appropriate temporal
   integration scheme, with common choices for FR being low storage explicit
   Runge-Kutta (RK) schemes~\cite{Kennedy2000,Vincent2011}. As was outlined in the introduction, the main concern of this paper is with the definition of the correction functions $h_L$ and $h_R$. The first few correction functions of FR were introduced by Huynh~\cite{Huynh2007}, which were broadly encapsulated into what we will call the Original stable Flux Reconstruction scheme (OSFR). This gave the definition of the correction functions as:
	\begin{align}
		h_L &= \frac{(-1)^p}{2}\bigg[ \psi_p - \Big(\frac{\eta_p\psi_{p-1} + \psi_{p+1}}{1+ \eta_p}\Big)\bigg] \label{eq:gl}\\
		h_R &= \frac{1}{2}\bigg[ \psi_p + \Big(\frac{\eta_p\psi_{p-1} + \psi_{p+1}}{1+ \eta_p}\Big)\bigg] \label{eq:gr}\\
	\end{align}
	where
	\begin{align}
		\eta_p &= \frac{\iota(2p+1)(a_pp!)^2}{2} \\
		a_p &= \frac{(2p)!}{2^p(p!)^2} 
	\end{align} 
	and $\psi_i$ is the $i^{\mathrm{th}}$ order Legendre polynomials of the first kind defined on $\xi\in[-1,1]$. Here the character of the scheme is defined by the free parameter $\iota$. The origin of this set of correction functions is from establishing that the scheme must be time stable in the Sobolev norm $|u + \iota(u^{(p)})|_{2}$, which is sufficient to define the topology of stable solutions. 
		
	An extended range of FR correction functions where then defined by Vincent~\etal~\cite{Vincent2015}, dubbed here ESFR. The derivation of this family finds it roots in DG and SE methods and the full proof is excluded for brevity. However, the key points are included here and we begin by defing the correction function and its gradient using a Legendre basis of the form:
	\begin{equation}\label{eq:dcorr_leg}
		h_L(\xi) = \sum^{p+1}_{i=0}\tilde{\mathbf{h_L}}_i\psi_i(\xi) \quad \mathrm{and} \quad \frac{dh_L(\xi)}{d\xi} = g_L(\xi) = \sum^p_{i=0} \tilde{\mathbf{g_L}}_i\psi_i(\xi)
	\end{equation}
	with the right correction similarly defined. Then the extended range of correction function is given by the equations:
	\begin{align}
		\tilde{\mathbf{g_L}} &= -\big(\tilde{\mathbf{M}} + \tilde{\mathbf{K}}\big)^{-1}\tilde{\mathbf{l}} \\
		\tilde{\mathbf{g_R}} &= \big(\tilde{\mathbf{M}} + \tilde{\mathbf{K}}\big)^{-1}\tilde{\mathbf{r}}
	\end{align} 
	where $\tilde{\mathbf{l}} = [\psi_0(-1) \dots \psi_p(-1)]^T$ and $\tilde{\mathbf{r}} = [\psi_0(1) \dots \psi_p(1)]^T$ give the values of the Legendre basis at the left and right interfaces. The mass matrix is defined as $\tilde{\mathbf{M}}$, given by:
	\begin{equation}
		\tilde{\mathbf{M}}_{i,j} = \int_{-1}^1 \psi_i\psi_j d\xi = \frac{2}{2j+1}\delta_{i,j}
	\end{equation} 
	$\tilde{\mathbf{K}}$ is defined to be our free variable, which must obey the following conditions:
	\begin{align}
		\tilde{\mathbf{K}} &= \tilde{\mathbf{K}}^{T} \\
		\tilde{\mathbf{K}}\tilde{\mathbf{D}} + (\tilde{\mathbf{K}}\tilde{\mathbf{D}})^T &= \mathbf{0} \\
		\tilde{\mathbf{M}} + \tilde{\mathbf{K}} &> \mathbf{0}
	\end{align}
	Together with $h_L(-1) = h_R(1) = 1$ and $h_L(1) = h_R(-1) = 0$, given the definitions:
	\begin{align}
		\mathbf{D}_{i,j} &= \frac{dl_j(\xi_i)}{d\xi} \\
		\mathbf{V}_{i,j} &= \psi_j(\xi_i) \\
		\tilde{\mathbf{D}} &= \mathbf{V}^{-1}\mathbf{D}\mathbf{V}
	\end{align}
	
	Therefore, this defines a multi-parameter family of correction functions which are enumerated and exemplified for various orders in Vincent~\etal~\cite{Vincent2015}. Lastly, OSFR and ESFR correction functions where grouped together by Trojak~\cite{Trojak2018} in the correction function set call Generalised Sobolev Stable FR, GSFR. The derivation of this set was made by exploring the functional Sobolev space more generally by consideration of all order of derivative up to the $p^{\mathrm{th}}$ derivative. For brevity the full stability proof is excluded, however GSFR can be adequately expressed as a left correction function that satisfies:
   \begin{equation}
       \sum^p_{i=0}\iota_i\int^1_{-1}\dxi{i}{h_L}{\xi}\pxi{i+1}{\hud}{\xi}d\xi
       = \sum^{p}_{i=1}\iota_i\Bigg|\pxi{i}{\hud}{\xi}\dxi{i}{h_L}{\xi}\Bigg|^1_{-1}
   \end{equation}
   where $\iota_i$ are $p+1$ free parameters that control the shape of the correction functions. For a full description of the means by which these correction functions are produced see Trojak~\cite{Trojak2018}, where also a matrix form is presented for easy evaluation from a given $\iota_i$.

\section{Generalised Lebesgue Stability}
	In the previous derivation of energy stable flux reconstruction, both OSFR and ESFR of Vincent~\etal~\cite{Vincent2010,Vincent2015} and GSFR of Trojak~\cite{Trojak2018}, the scheme was shown to be stable in the broken modified Sobolev norm. This results in a scheme whose energy is stable when balanced between the conserved variable and derivatives of the conserved variable. However, Trojak~\cite{Trojak2018} notes that in the case when $\iota_i=0 \:\forall\: i\ne0$ there is an illposedness in the functional space defined by GSFR. Hence, we wish to explore the stability of this case by setting the Sobolev norm appropriately, which results in the familiar Lebesgue norm. First introducing the modified broken Sobolev norm:
	\begin{equation}
		\|u\|_{n,W^{p,2,\iota}} = \sqrt{\int_{\Omega_n} \sum^p_{i=0}\iota_i\bigg(\dxi{i}{u}{\xi}\bigg)^2 d\xi}
	\end{equation}
	which, when we take $\iota_i=0 \:\forall\: i\ne 0$ and $\iota_0 = 1$, this becomes the broken Lebesgue norm:
	\begin{equation}
		\|u\|_{n,L^2} = \sqrt{\int_{\Omega_n} u^2 d\xi}
	\end{equation}
	We will now look for corrections that satisfy the stability criterion, $d(\|u\|_{n,L^2})/dt \leqslant 0$, in this norm for the 1D FR conservation law:
	\begin{equation}
		\px{\hud}{t} = -\px{\hfd}{\xi} - (\hfI_L - \hfd_L)\dx{h_L}{\xi} - (\hfI_R - \hfd_R)\dx{h_R}{\xi}
	\end{equation}
	If the flux function is set such that linear advection is solved then $\hfd = \hud$. This then implies that:
	\begin{equation} \label{eq:fr_linear}
      \px{\hud}{t} = -\px{\hud}{\xi} - (\huI_L - \hud_L)\dx{h_L}{\xi} - (\huI_R - \hud_R)\dx{h_R}{\xi}
	\end{equation}
	Multiplying Eq.(\ref{eq:fr_linear}) by the discontinuous conserved variable and integrating over the sub-domain:
	\begin{equation} \label{eq:fr_strong}
		\int^1_{-1}\hud\px{\hud}{t} d\xi = -\int^1_{-1}\hud\px{\hud}{\xi}d\xi - (\huI_L - \hud_L)\int^1_{-1}\hud\dx{h_L}{\xi}d\xi - (\huI_R - \hud_R)\int^1_{-1}\hud\dx{h_R}{\xi}d\xi
	\end{equation}
	Applying the product rule and integration by parts to Eq.(\ref{eq:fr_strong}), the rate of energy decay of FR which may also be recognised as the rate of change of the Lebesgue norm is:
	\begin{equation}\label{eq:FR_0sobolev}
	 	\frac{1}{2}\dx{}{t}\int^1_{-1}(\hud)^2 d\xi = -\frac{1}{2}\int^1_{-1}\px{(\hud)^2}{\xi} d\xi - \big(\huI_L - \hud_L\big)\underbrace{\int^1_{-1} h_L\px{\hud}{\xi}d\xi}_{I_L} - \big(\huI_R-\hud_R\big)\underbrace{\int^1_{-1} h_R\px{\hud}{\xi} d\xi}_{I_R} + \big(\huI_L\hud_L - \huI_R\hud_R\big)
	\end{equation}
	If the integrals $I_L$ and $I_R$ are set equal to zero, then Eq.(\ref{eq:FR_0sobolev}) may be reduced to:
	\begin{equation}
		\frac{1}{2}\dx{}{t}\int^1_{-1}(\hud)^2 d\xi = -\frac{1}{2}\big((\hud_L)^2 - (\hud_R)^2\big) + \big(\huI_L\hud_L - \huI_R\hud_R\big)
	\end{equation}
	This can be further simplified into an expression for the transfer of energy across the boundary.
	\begin{equation}
		\frac{1}{2}\dx{}{t}\int^1_{-1}(\hud)^2 d\xi = \frac{1}{2}\Big(2\huI_L-\hud_L\Big)\hud_L - \frac{1}{2}\Big(2\huI_R-\hud_R\Big)\hud_R
	\end{equation}
	and hence the energy stability after setting of the correction is controlled by the method used for interface calculation. In order to extract a condition in which $I_L=I_R=0$, it is useful to consider the orthogonal polynomial basis, by taking $\hud$ to be:
	\begin{equation}\label{eq:leg_u}
		\hud(\xi;t) = \sum^p_{i=0}\hat{v}_i(t)\psi_i(\xi)
	\end{equation}
	and using the similar basis for the correction function as in Eq.(\ref{eq:dcorr_leg}). Is may be useful at this point to be reminded of the result for Legendre polynomials~\cite{Holdeman1970}, caused by there alternating odd-even nature:
	\begin{equation} \label{eq:legendre_udu}
		\int^1_{-1} \psi_l(\xi)\dx{\psi_m(\xi)}{\xi} d\xi = \begin{cases}
2, \quad \mathrm{for} \quad l/2 - \floor{l/2} \neq m/2 - \floor{m/2}, \quad l\leqslant m\\
0, \quad \mathrm{otherwise} \\ 
		\end{cases}
	\end{equation}
	Hence, applying this to $I_L$ and $I_R$, the family of correction functions is found by satisfaction of the following conditions:
	\begin{equation}\label{eq:GLSFR_condition}
		\sum_{\substack{i=0\\ i=\mathrm{even}}}^{p-1} \tilde{\mathbf{h_l}}_i = 0 \quad \mathrm{and} \quad
		\sum_{\substack{i=0\\ i=\mathrm{odd}}}^{p-1} \tilde{\mathbf{h_l}}_i = 0 \\
	\end{equation}
	By enforcing the boundary conditions namely: $h_L(-1) = 1$ and $h_L(1) = 0$, the final conditions on $\tilde{\mathbf{h_L}}$ may be found as:
	\begin{equation}
		\tilde{\mathbf{h_l}}_{p+1} = \frac{(-1)^{p+1}}{2} \quad \mathrm{and} \quad \tilde{\mathbf{h_l}}_p = \frac{(-1)^{p}}{2}
	\end{equation}
	It can be seen by comparison to the OSFR scheme of Eq.(\ref{eq:gl}), that this new set of correction functions are coincident with the previous family of corrections at only one point, $\eta_p=0$, corresponding to Nodal DG or when the terms of Eq.(\ref{eq:GLSFR_condition}) are all zero. Owing to space of stability this scheme inhabits, this new set of corrections will henceforth be called Generalised Lebesgue Stable FR, GLSFR.
		
\section{Conservation of Mass}
	It may at this point be thought that the highly arbitrary family of correction functions that satisfy Eq.\ref{eq:GLSFR_condition} may lead to the injection or sink of mass into the sub-elements. To placate this fear we consider the integration of Eq.(\ref{eq:FR_broken}) over the sub-domain, \emph{i.e.} the rate of change of mass in the sub-domain:
	\begin{equation}
		\dx{}{t}\int^{1}_{-1}\hud d\xi = -\int^{1}_{-1}\px{\hfd}{\xi} d\xi - \Big(\hfI_L-\hfd_L\Big)\int^1_{-1}\dx{h_L}{\xi}d\xi - \Big(\hfI_R-\hfd_R\Big)\int^1_{-1}\dx{h_R}{\xi} d\xi
	\end{equation} 
	and considering the constraints on $h_L$ and $h_R$, this becomes:
	\begin{align}
		\dx{}{t}\int^{1}_{-1}\hud d\xi &= \Big(\hfd_L-\hfd_R\Big) + \Big(\hfI_L-\hfd_L\Big) - \Big(\hfI_R-\hfd_R\Big) \\
		&= \hfI_L - \hfI_R
	\end{align}
		Hence, the rate of mass accumulation within an element is dependant on the interface calculation, the flux function and, as was shown in Trojak~\etal~\cite{Trojak2017}, mesh deformation. But is not dependant on the correction function, so long as the boundary conditions are fully enforced.

\section{von Neumann Analysis}
	Further investigation into the stability and performance of this set of correction functions is to be performed via von Neumann analysis. Where we will inspect the spectral properties and ergodicity of FR applied to the linear advection and advection-dffusion equations for harmonic solutions. This form of analysis was outlined by Lele~\cite{Lele1992}, for finite difference, applied to DG by Hu~\etal~\cite{Hu1999,Hesthaven2008} and adapted to FR by Huynh~\cite{Huynh2007} and Vincent~\etal~\cite{Vincent2011}. The formulation and structure of the analysis presented here follows that of Trojak~\etal~\cite{Trojak2017} and we may start by defining the semi-discrete matrixised form of FR applied to the linear advection equation:
	\begin{equation}\label{eq:la_fr_matrix}
		\px{\mathbf{u}_j}{t} = - J_{j+1}^{-1}\mathbf{C}_{+}\mathbf{u}_{j+1} - J_{j}^{-1}\mathbf{C}_0\mathbf{u}_j - J_{j-1}^{-1}\mathbf{C}_{-}\mathbf{u}_{j-1}
	\end{equation}
	where the matrixised operators are defined as:
	\begin{align}
		\mathbf{C}_{+} &= (1-\alpha)\mathbf{g_rl_l}^T \\
		\mathbf{C}_0  &= \mathbf{D} - \alpha\mathbf{g_ll_l}^T - (1-\alpha)\mathbf{g_rl_r}^T \\
		\mathbf{C}_{-} &= \alpha\mathbf{g_ll_r}^T
	\end{align}
	where $\mathbf{D}$ is the discrete differentiation operator, $\mathbf{g_l}$ and $\mathbf{g_r}$ are respectively the left are right correction function gradients at the solution points, $\mathbf{l_l}$ and  $\mathbf{l_r}$ are respectively the interpolation weights to the left and right flux points from the solution points. $\alpha$ is then defined to be the upwinding ratio, ($\alpha =1$ meaning fully upwinded and $\alpha = 0.5$ meaning central differenced). A harmonic solution can then be imposed on this by using the bloch wave:
	\begin{equation}\label{eq:bloch}
		u = v\exp{\big(i(kx-\omega t)\big)}
	\end{equation}
	which, upon substitution into Eq.(\ref{eq:la_fr_matrix}), the harmonic solution to  linear advection via FR is thus:
	\begin{equation}\label{eq:FR_ad_semi}
		\px{\mathbf{u}_j}{t}  = -\Big(J^{-1}_{j+1}\mathbf{C}_{+}\exp{(-ik\delta_j)} +J^{-1}_{j}\mathbf{C}_{0} + -J^{-1}_{j-1}\mathbf{C}_{-}\exp{(-ik\delta_{j-1})}\Big)\mathbf{u}_j = \mathbf{Q}_{\mathrm{a}}\mathbf{u}_j
	\end{equation}
	
	Taking this analysis one set further, most problem of practical interest involve second order derivatives. Therefore, we wish to understand the behaviour and stability of this branch of correction functions when applied to diffusion and advection-diffusion problems. If we introduce the linear diffusion equation written as:
	\begin{align}
		\px{u}{t} &= \nu\px{q}{x} \\
		q &= \px{u}{x} 
	\end{align}
	Thus from Eq.(\ref{eq:FR_ad_semi}) we may write:
	\begin{equation}
		\px{\mathbf{u}_j}{t} = J^{-1}_{j-1}\mathbf{C}_{-1}J^{-1}_{j-1}\mathbf{q}_{j-1}\: +\: J^{-1}_{j}\mathbf{C}_{0}J^{-1}_{j}\mathbf{q}_{j}\: +\: J^{-1}_{j+1}\mathbf{C}_{+1}J^{-1}_{j+1}\mathbf{q}_{j+1} 
	\end{equation}
	By applying the same spatial discretisation to $\mathbf{q}$, $\mathbf{q}_j$ can be written as:
	\begin{equation}
		J^{-1}_{j} \mathbf{q}_j = 
		J^{-1}_{j-1}\mathbf{C}_{-1}\mathbf{u}_{j-1} \: + \:
		J^{-1}_{j}\mathbf{C}_{0}\mathbf{u}_{j} \: + \:
		J^{-1}_{j+1}\mathbf{C}_{+1}\mathbf{u}_{j+1}
	\end{equation}
	This implies that the same correction function is used for both the diffusion correction and advection correction. This was found to give optimal performance~\cite{Williams2014a,Castonguay2012} and has the benefit of easier practical implementation.Proceeding, the semi-discretised linear diffusion equation for FR is then:
	\begin{equation}\label{eq:FR_D2_semi_disc}
		\begin{split}
		\px{\mathbf{u}_j}{t} =&\:\:
		\underbrace{J^{-1}_{j-2}J^{-1}_{j-1}\mathbf{C}_{-1}^2}_{\mathbf{B}_{-2}}\mathbf{u}_{j-2}\:\: + \\
		&\underbrace{J^{-1}_{j-1}(J^{-1}_{j-1}\mathbf{C}_{-1}\mathbf{C}_{0}\:\: + \:\:J^{-1}_{j}\mathbf{C}_{0}\mathbf{C}_{-1})}_{\mathbf{B}_{-1}}\mathbf{u}_{j-1}\:\: + \\
		&\underbrace{J^{-1}_{j}(J^{-1}_{j-1}\mathbf{C}_{-1}\mathbf{C}_{+1}\:\: + \:\:J^{-1}_{j}\mathbf{C}_{0}^2\:\: + \:\:J^{-1}_{j+1}\mathbf{C}_{+1}\mathbf{C}_{-1})}_{\mathbf{B}_{0}}\mathbf{u}_{j}\:\: + \\
		&\underbrace{J^{-1}_{j+1}(J^{-1}_{j}\mathbf{C}_{0}\mathbf{C}_{+1}\:\: + \:\:J^{-1}_{j+1}\mathbf{C}_{+1}\mathbf{C}_{0})}_{\mathbf{B}_{+1}}\mathbf{u}_{j+1}\:\: + \\
		&\underbrace{J^{-1}_{j+1}J^{-1}_{j+2}\mathbf{C}_{+1}^2}_{\mathbf{B}_{+2}}\mathbf{u}_{j+2}
		\end{split}
	\end{equation}
	where the matrices $\mathbf{C}_{-1}$, $\mathbf{C}_0$ and $\mathbf{C}_{+1}$ are defined as before. This form is very similar to that shown in \cite{Watkins2016}. In this case, we will use BR1 for the  calculation of the common interface flux, hence $\alpha=0.5$. Applying a trial solution of $u = v\exp{(ikx_j)}\exp{(-k^2c_dt)}$ to Eq.(\ref{eq:FR_D2_semi_disc}) and simplifying:
	\begin{equation}
		c_d\mathbf{v} =-\frac{1}{k^2}\Big(e^{ik(x_{j-2}-x_j)}\mathbf{B}_{-2}\mathbf{v} + e^{ik(x_{j-1}-x_j)}\mathbf{B}_{-1}\mathbf{v} + \mathbf{B}_{0}\mathbf{v} + e^{ik(x_{j+1}-x_j)}\mathbf{B}_{+1}\mathbf{v} + e^{ik(x_{j+2}-x_j)}\mathbf{B}_{+2}\mathbf{v}\Big)
	\end{equation}	
	If $\delta_j = x_j-x_{j-1}$, then this can be further simplified to:
	\begin{equation}\label{eq:FR_D2_VN}
		c_d\mathbf{v} =-\frac{1}{k^2}\Big(e^{-ik(\delta_{j-1}+\delta_j)}\mathbf{B}_{-2} + e^{-ik\delta_j}\mathbf{B}_{-1} + \mathbf{B}_{0} + e^{ik\delta_{j+1}}\mathbf{B}_{+1} + e^{ik(\delta_{j+2}+\delta_{j+1})}\mathbf{B}_{+2}\Big)\mathbf{v}
	\end{equation}
	Which can also allow us to define the FR diffusion matrix as:
	\begin{equation}\label{eq:FR_D2_Q}
		\mathbf{Q}_{\mathrm{d}} = \Big(
		e^{-ik(\delta_{j-1}+\delta_j)}\mathbf{B}_{-2} +	e^{-ik\delta_j}\mathbf{B}_{-1} + \mathbf{B}_{0}  + e^{ik\delta_{j+1}}\mathbf{B}_{+1} + 
		e^{ik(\delta_{j+2}+\delta_{j+1})}\mathbf{B}_{+2}\Big)
	\end{equation}
	
	Equation.~(\ref{eq:FR_D2_VN}) is again an eigenvalue problem, albeit a non-trivial one, where from \cite{Trefethen1994,Watkins2016} the dissipation and dispersion are defined as $\Re{(\hat{k}^2c_d)}$ and  $\Im{(\hat{k}^2c_d)}$ respectively, due to the  of the second derivative. The physical mode is extracted from the $p+1$ dimensional eigenproblem following the procedure of \cite{Watkins2016}. Finally for the case of linear advection-diffusion, the fully discrete form may be written using the update matrix, which encompasses the temporal integration method.
	\begin{align}\label{eq:update}
		\mathbf{Q}_{\mathrm{ad}} &= 2c\mathbf{Q}_{\mathrm{a}} + 4\nu\mathbf{Q}_{\mathrm{d}} \\
		\mathbf{u}^{n+1}_j &= \mathbf{R}(\mathbf{Q}_{\mathrm{ad}})\mathbf{u}^n_j \\
		\mathbf{R}_{33}(\mathbf{Q}_{\mathrm{ad}}) &= \sum_{m=0}^3 \frac{(\tau\mathbf{Q}_{\mathrm{ad}})^m}{m!}	
	\end{align}
	From the the update matrix $\mathbf{R}$ the fully discrete dispersion and dissipation can then be found but further substitution of the trial solution of Eq.(\ref{eq:bloch}) into Eq.(\ref{eq:update}).
	\begin{align}
		\underbrace{e^{-ik(c-1)\tau}}_{\lambda}\mathbf{v} &= e^{ik\tau}\mathbf{Rv}\label{eq:full_discrete} \\
		c(k;\tau) &= \frac{i\log{(\lambda)}}{k\tau} + 1 
	\end{align}
	Hence, from von Neumann's theorem~\cite{Isaacson1994}, for temporal stability of the fully discrete scheme the spectral radius of the update matrix must be less than of equal to one, $\rho(\mathbf{R}) \leqslant 1$. For comparative purposes we then define the normalised time step for advection-diffusion as:
	\begin{equation}
		\hat{\tau} = \bigg(\frac{2c}{h} + \frac{4\nu}{h^2}\bigg)\tau
	\end{equation}
   and this will be used to define the CFL limit.

\subsection{Results}
	Initially we will study the dispersion and dissipation of linear advection, with the primary aim of showing how, for what may be traditionally considered an inappropriate correction function, we are able to recover a stable scheme. Take the following example, when $p=4$ and the correction weights are chosen arbitrarily to be:
	\begin{align}
		\tilde{\mathbf{h_l}}_0 &= (5.22943203125\times10^4)\times10^{-5} \label{eq:h0_test1}\\
		\tilde{\mathbf{h_l}}_1 &= 0.1\sqrt{2} \label{eq:h1_test1}
	\end{align}
	where Eq.(\ref{eq:h0_test1}) is $10^{-5}$ times the floating point representation of the ASCII string `\texttt{GLFR}`.
	\begin{figure}[h]
	 	\centering
		\begin{subfigure}[b]{0.32\linewidth}
			\centering
			\includegraphics[width=\linewidth]{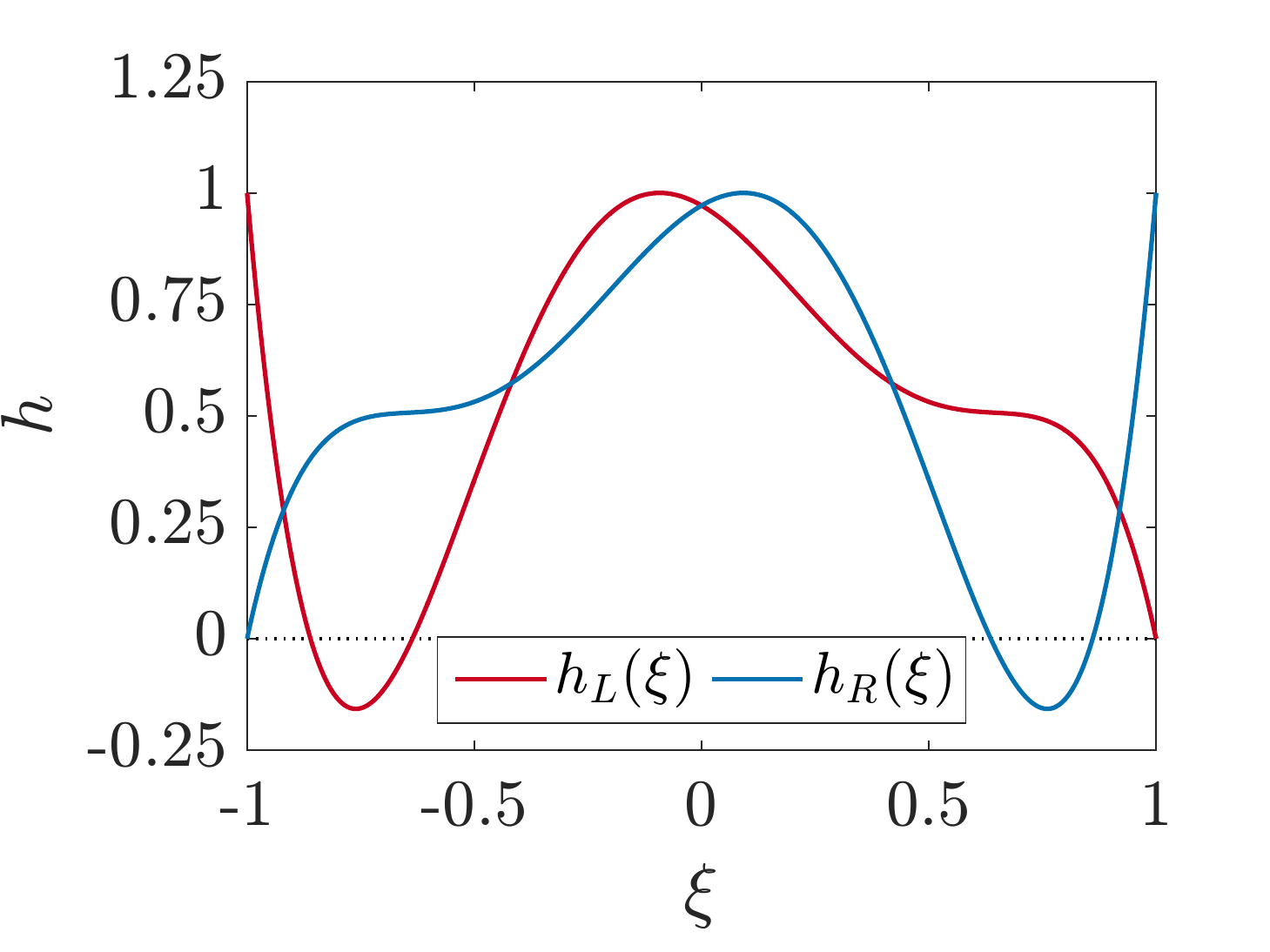}
			\caption{Left and right correction functions}
			\label{fig:FR4GLS_random_corr}
		\end{subfigure}
		~
		\begin{subfigure}[b]{0.32\linewidth}
			\centering
			\includegraphics[width=\linewidth]{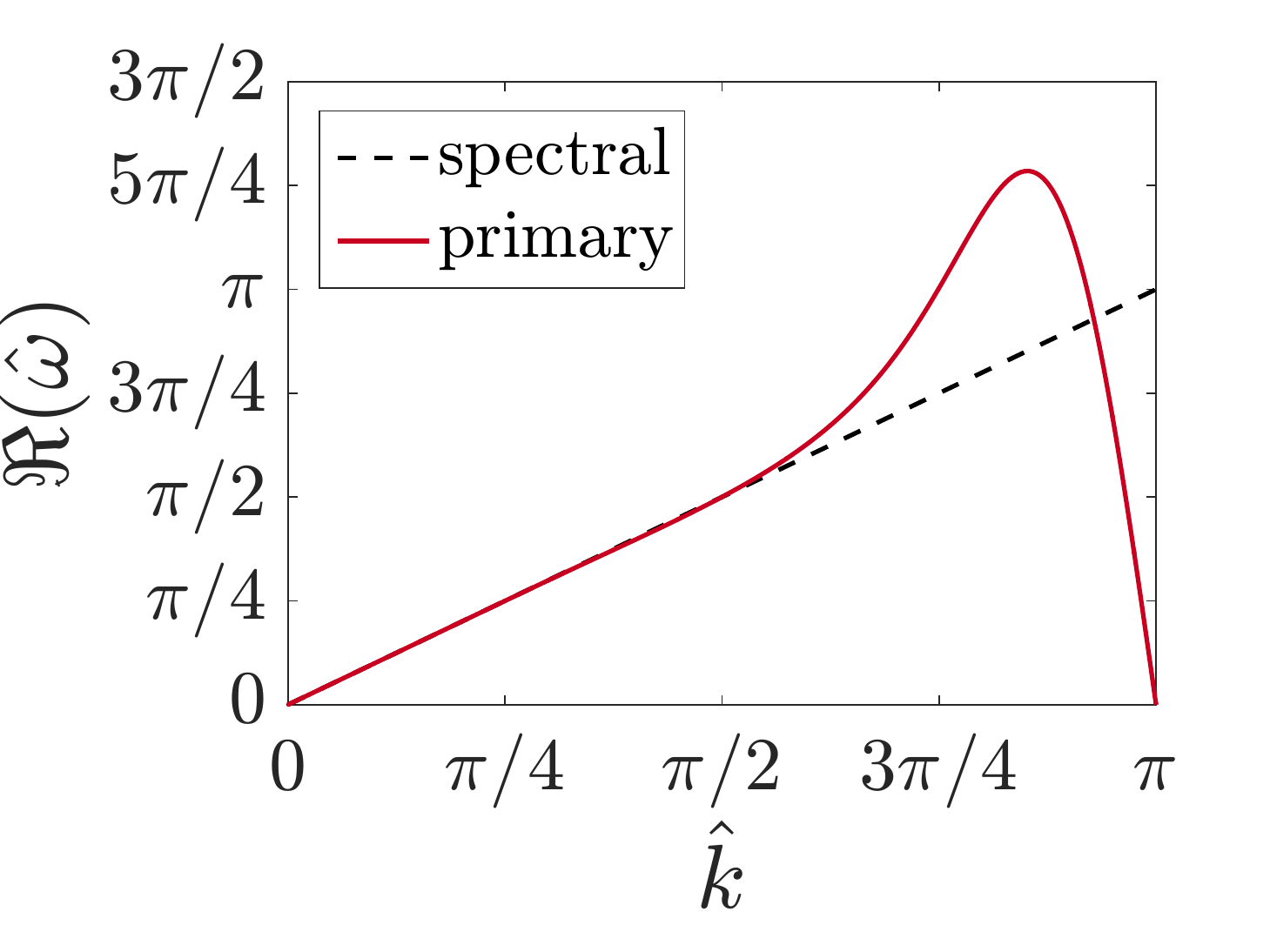}
			\caption{dispersion}
			\label{fig:FR4GLS_random_R}
		\end{subfigure}
		~
		\begin{subfigure}[b]{0.32\linewidth}
			\centering
			\includegraphics[width=\linewidth]{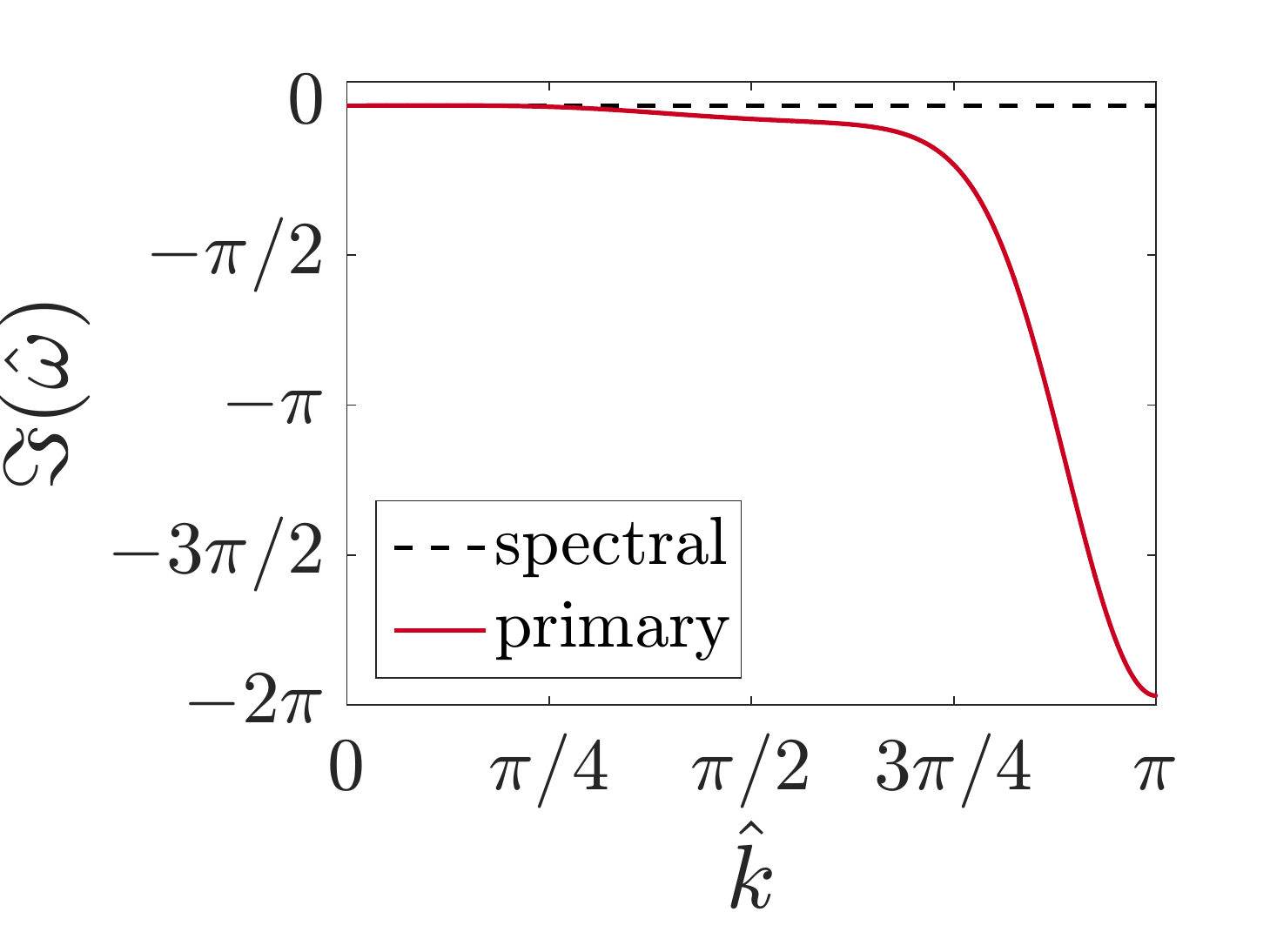}
			\caption{dissipation}
			\label{fig:FR4GLS_random_I}
		\end{subfigure}
		\caption{Correction function and wave propagating characteristics for arbitrary correction weights as in Eq.(\ref{eq:h0_test1}~\&~\ref{eq:h1_test1}), for $p=4$ on a uniform grid with upwinded cell interfaces applied to linear advection.}
		\label{fig:FR4GLS_random}
	\end{figure}
	It is clearly apparent that this arbitrary correction function that satisfies Eq.~(\ref{eq:GLSFR_condition}) is stable and although wave propagating characteristic of this arbitrary example are not optimal. However, this aims to demonstrate initially that a correction function may be defined that was previously unproducible, leading to an incredibly general family of schemes.
	
	In order to understand the space of temporally stable GLSFR correction functions the 1D linear-advection von Neumann analysis is first investigated for various orders, correction functions and interface calculations.
	\begin{figure}
		\centering
		\begin{subfigure}[b]{0.45\linewidth}
			\centering
			\includegraphics[width=0.9\linewidth]{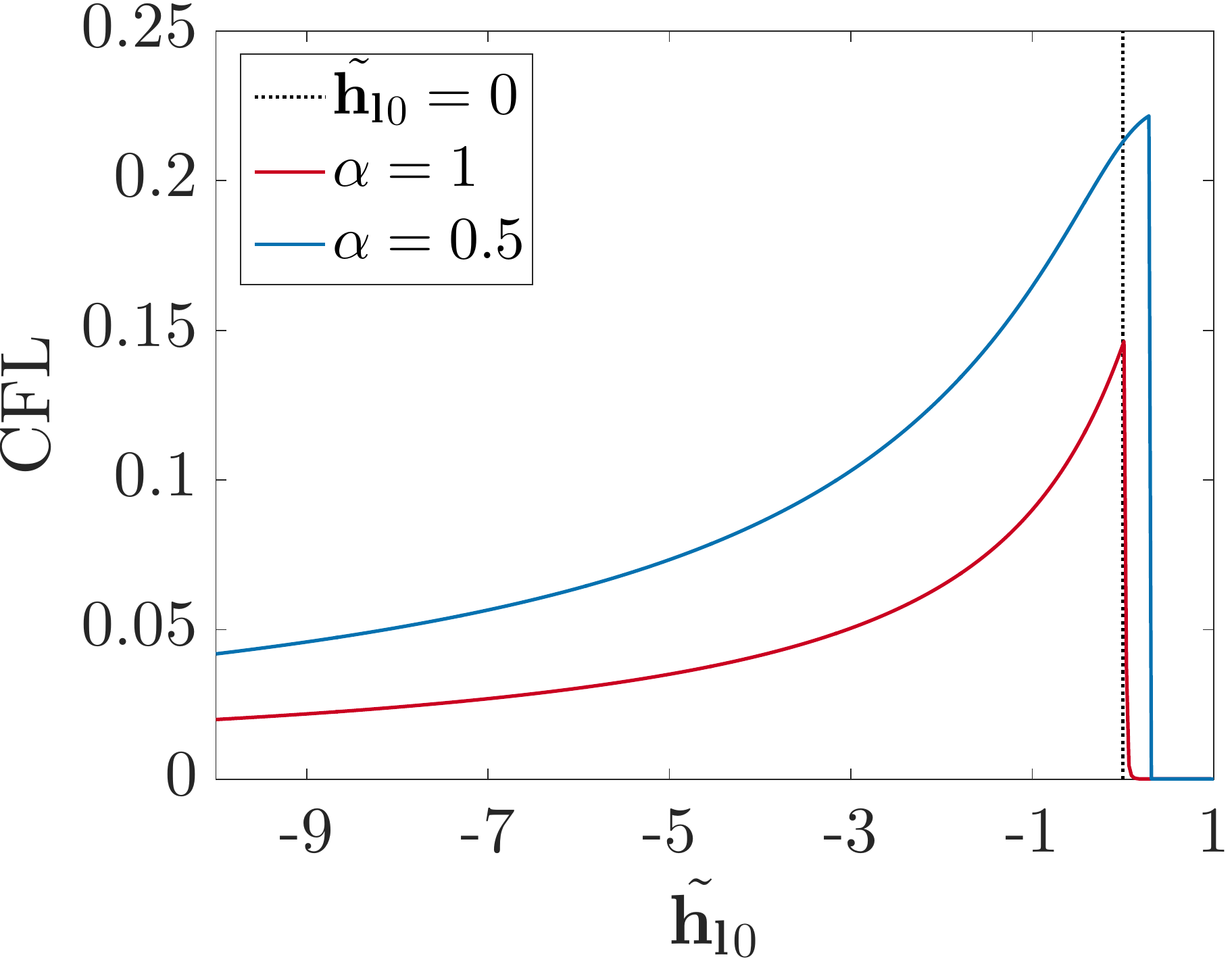}
			\caption{$p=3$}
			\label{fig:FRGLS3_RK44_CFL}
		\end{subfigure}
		~\\
		\begin{subfigure}[b]{0.45\linewidth}
			\centering
			\includegraphics[width=\linewidth]{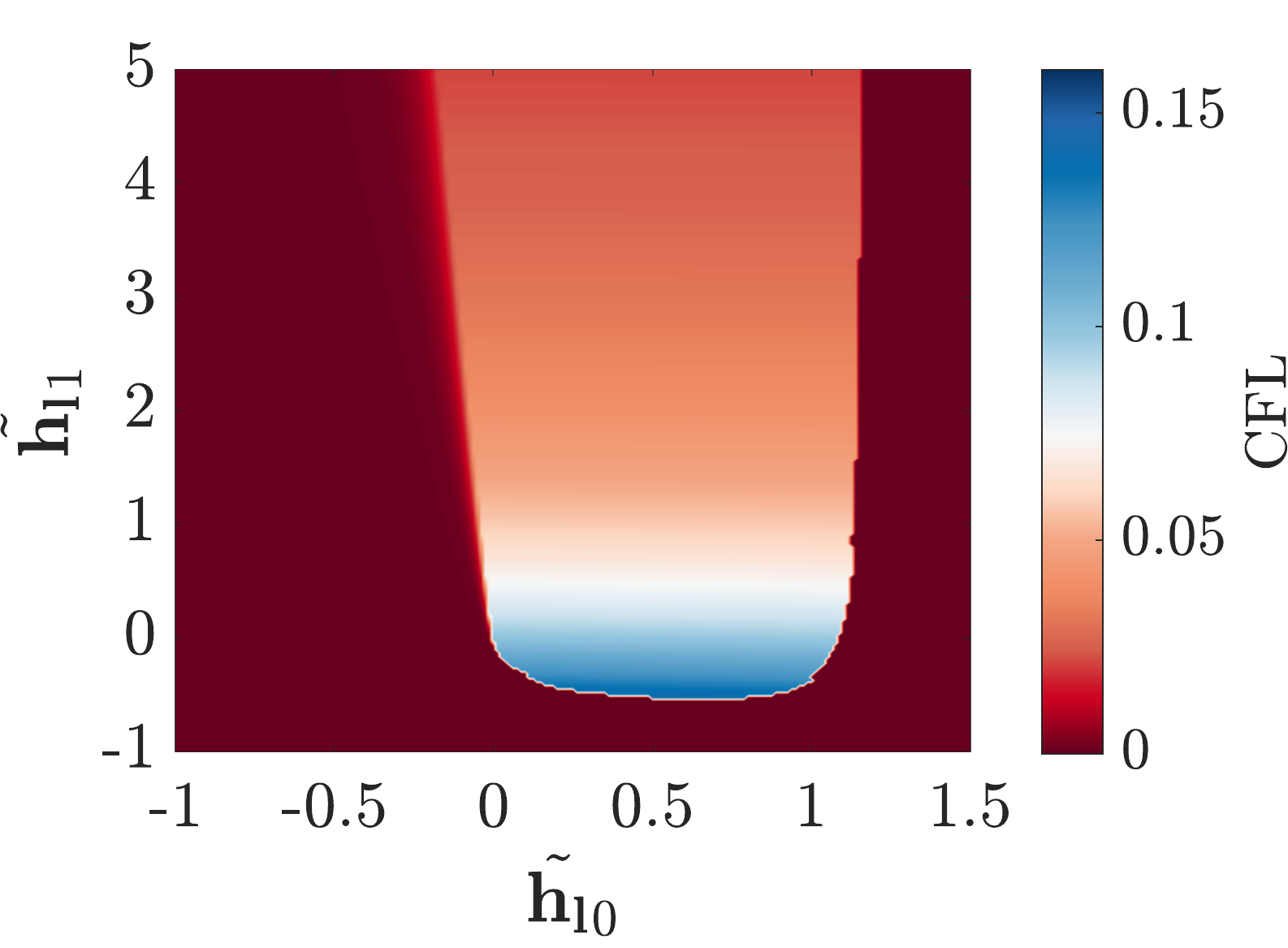}
			\caption{$p=4$, $\alpha = 1$ (upwinded)}
			\label{fig:FRGLS4_RK44_CFL}
		\end{subfigure}
		~
		\begin{subfigure}[b]{0.45\linewidth}
			\centering
			\includegraphics[width=\linewidth]{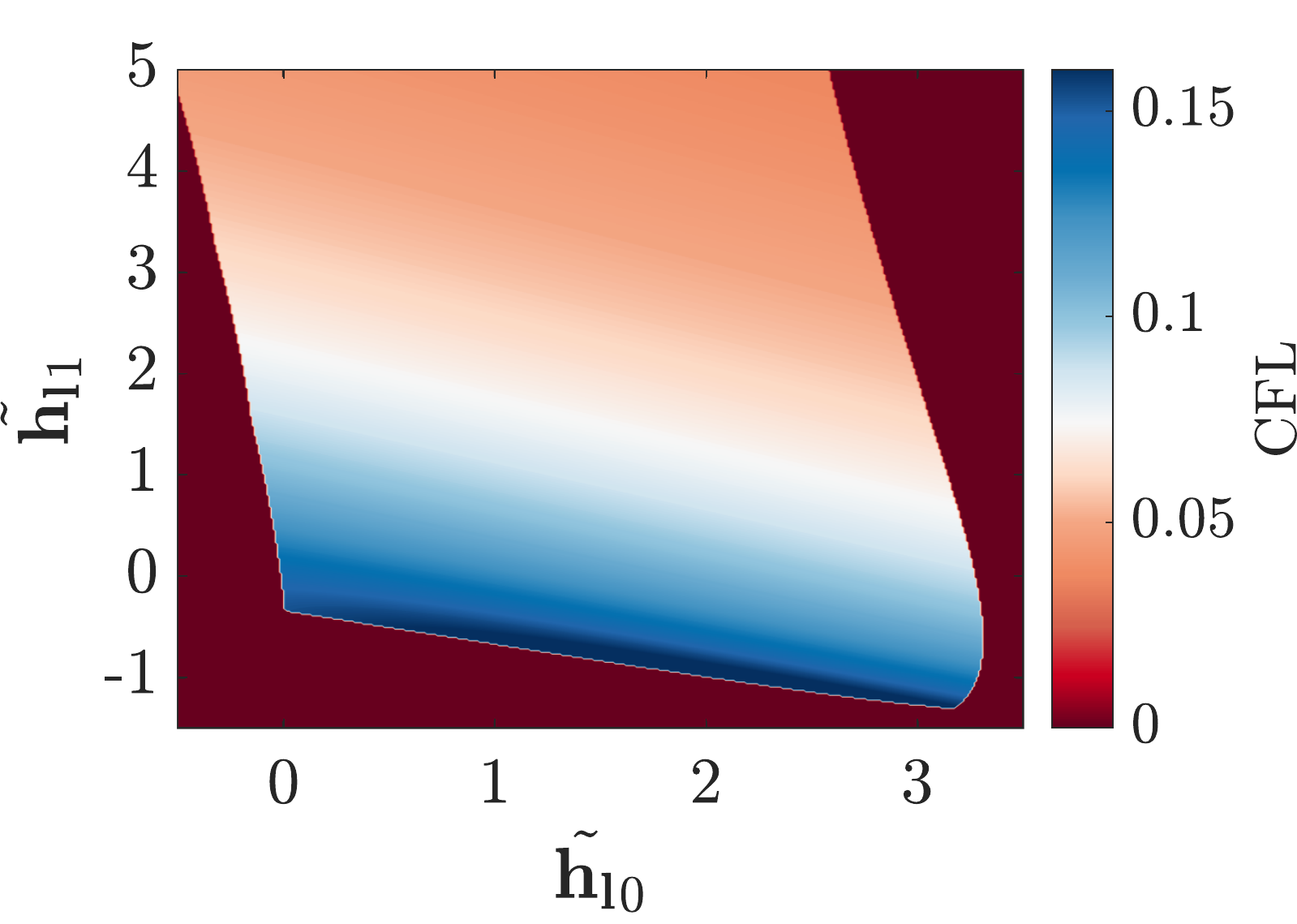}
			\caption{$p=4$, $\alpha = 0.5$ (central)}
			\label{fig:FRGLS4_RK44_CFL_central}
		\end{subfigure}
		\caption{Maximum stable CFL number for GLSFR $[c=1,\nu=0]$, with RK44 temporal integration.}
		\label{fig:GLSFR_CFL}
	\end{figure}
	
	Figure~\ref{fig:GLSFR_CFL} shows the region of stable correction function when RK44 temporal integration is used, with a clear region of overlap between the upwinded and central cases. For previous correction function sets, the existence of the Sobolev norm was used to determine the region of stable correction functions, however due to the nature of the GLSFR Sobolev norm, in that it collapses onto the $L^2$ norm, this approach could not be taken. Instead the results of Fig.\ref{fig:GLSFR_CFL} may be used to heuristically bound the correction function set.	To demonstrate the wave propagation properties of GLSFR we perform a semi-discrete and fully discrete von Neumann analysis. We will take the case when $p=4$ and compare the correction function that gives the highest CFL limit for both the upwinded and central case to DG. We chose DG as this is a commonly used scheme even within FR~\cite{Vincent2017} due to the complexity currently with the simplex correction definition~\cite{Castonguay2012a}. The CFl optimal case corresponds to $\mathbf{H}_4 =[\tilde{\mathbf{h_l}}_0=0.77,\tilde{\mathbf{h_l}}_1=-0.52]$, with the DG comparison at $\mathbf{H}_4=[0,0]$. 
	
	\begin{figure}[h]
	 	\centering
		\begin{subfigure}[b]{0.45\linewidth}
			\centering
			\includegraphics[width=\linewidth]{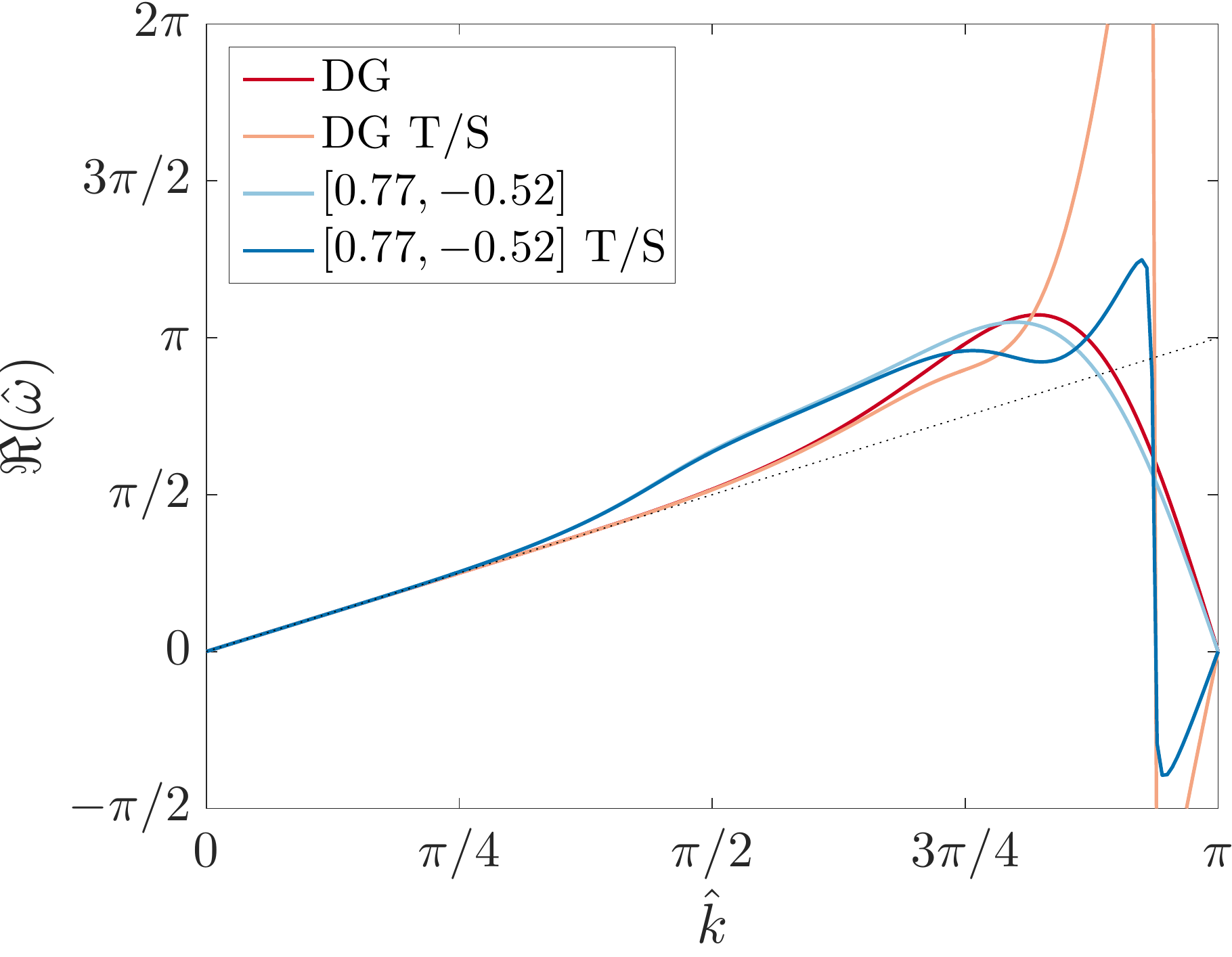}
			\caption{Dispersion}
			\label{fig:FRGLS4_dg_disp}
		\end{subfigure}
		~
		\begin{subfigure}[b]{0.45\linewidth}
			\centering
			\includegraphics[width=\linewidth]{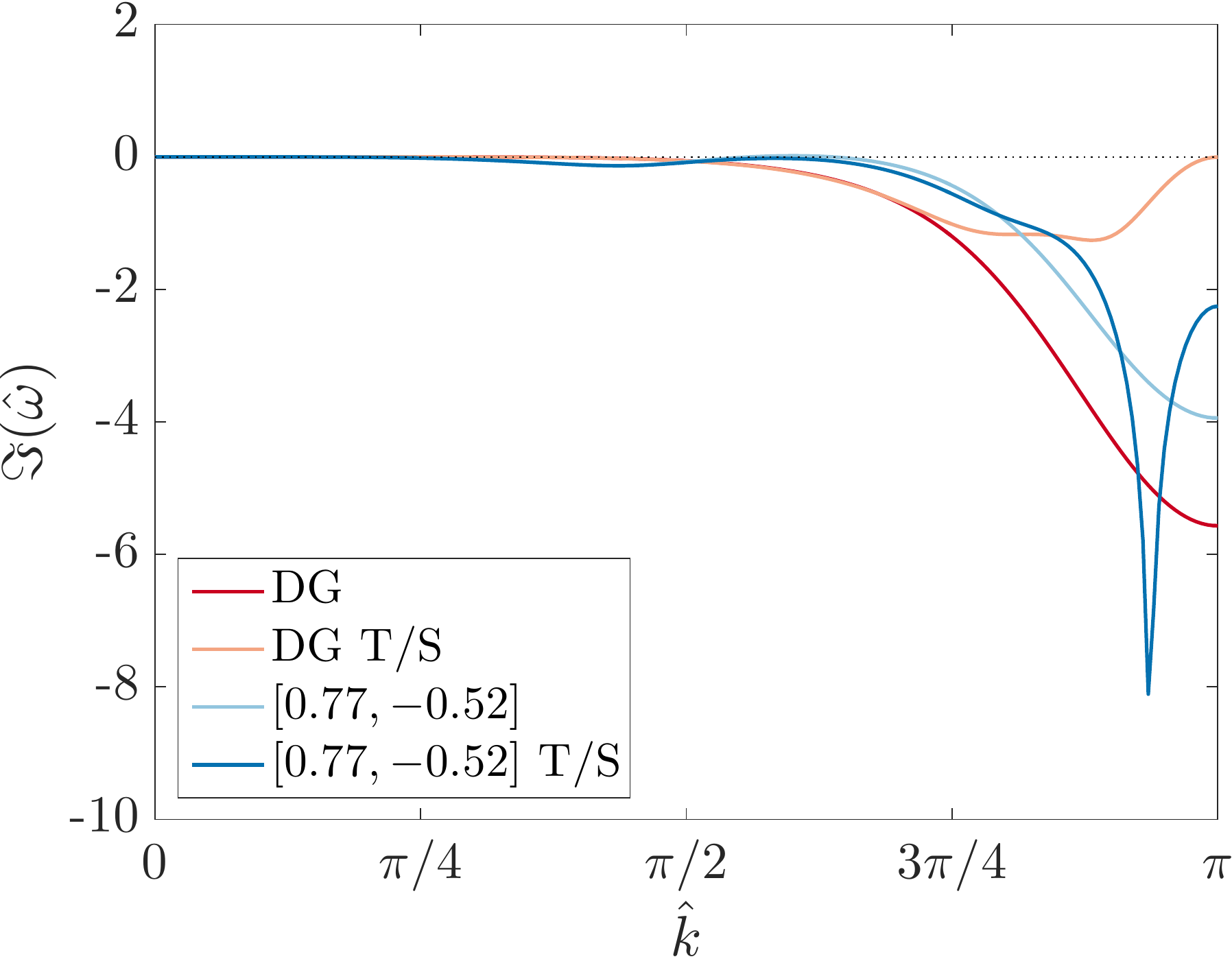}
			\caption{Dissipation}
			\label{fig:FRGLS4_dg_diss}
		\end{subfigure}
		\caption{Linear advection dispersion and dissipation comparison of spatial and temporal-spatial (T/S) analysis for upwinded FR, $p=4$. This compares DG and with GLSFR with $\mathbf{H}_4 = [0.77,-0.52]$. The temporal-spatial analysis uses RK44 temporal integration at $\tau = 0.1$ in both cases.}
		\label{fig:FRGLS_disp_diss_DG}
	\end{figure} 	
	
	It should be clear from Fig.\ref{fig:FRGLS_disp_diss_DG} that GLSFR in the fully discrete case is able to greatly reduce the wave group velocity, however at the expense of dispersion overshoot. Furthermore, by comparison of the dispersion and dissipation it can be seen that this GLSFR correction function has a very localised region of high dissipation around the wavenumber where the dispersion leads us to have a very large negative group velocity. It can be concluded that such a scheme may give improved performance when applied to explicit LES.
	
	\begin{figure}[h]
	 	\centering
		\begin{subfigure}[b]{0.45\linewidth}
			\centering
			\includegraphics[width=\linewidth]{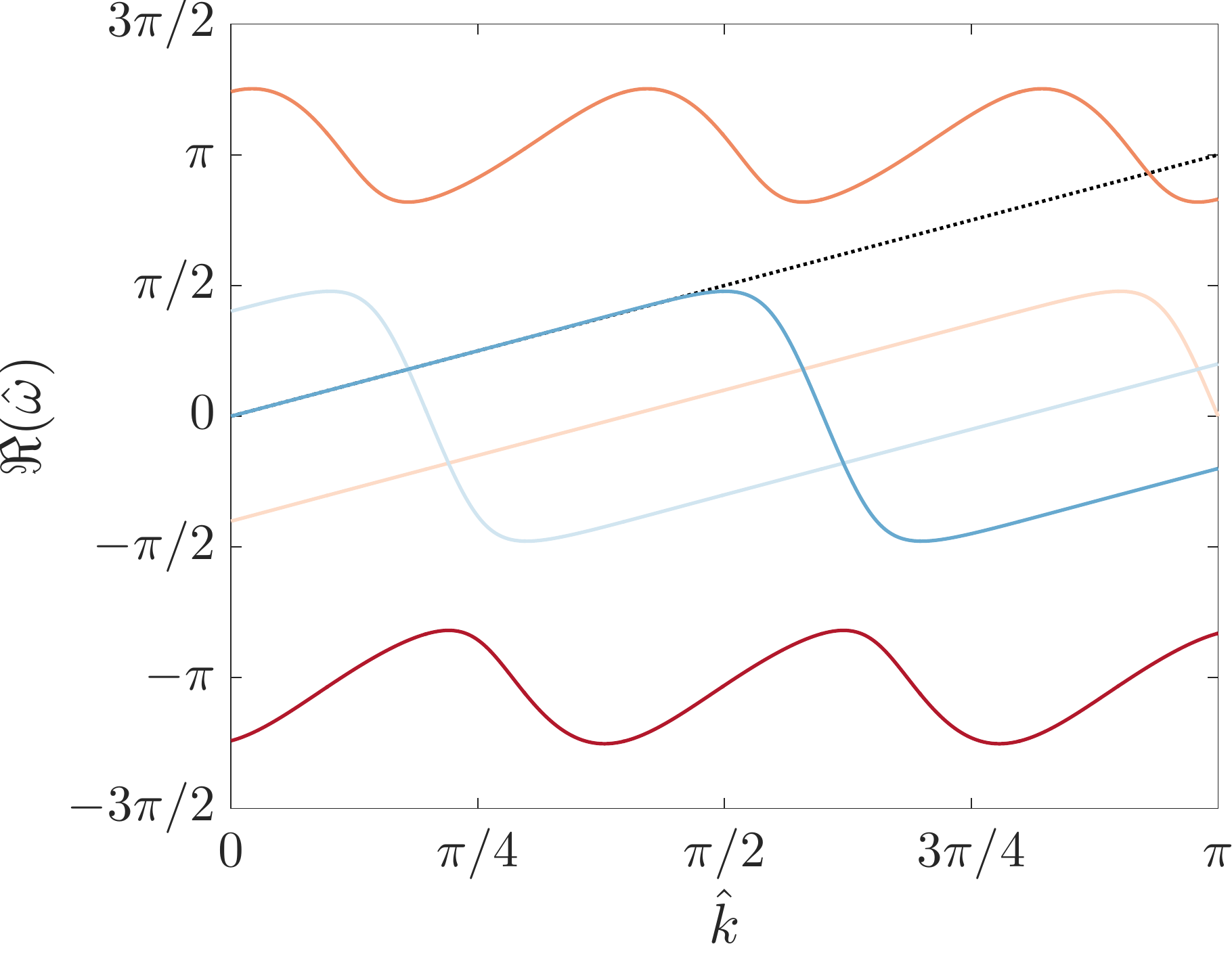}
			\caption{$\mathbf{H}_4 = [0,0]$, equivalent to DG}
			\label{fig:FRGLS4_dg_diss_c}
		\end{subfigure}
		~
		\begin{subfigure}[b]{0.45\linewidth}
			\centering
			\includegraphics[width=\linewidth]{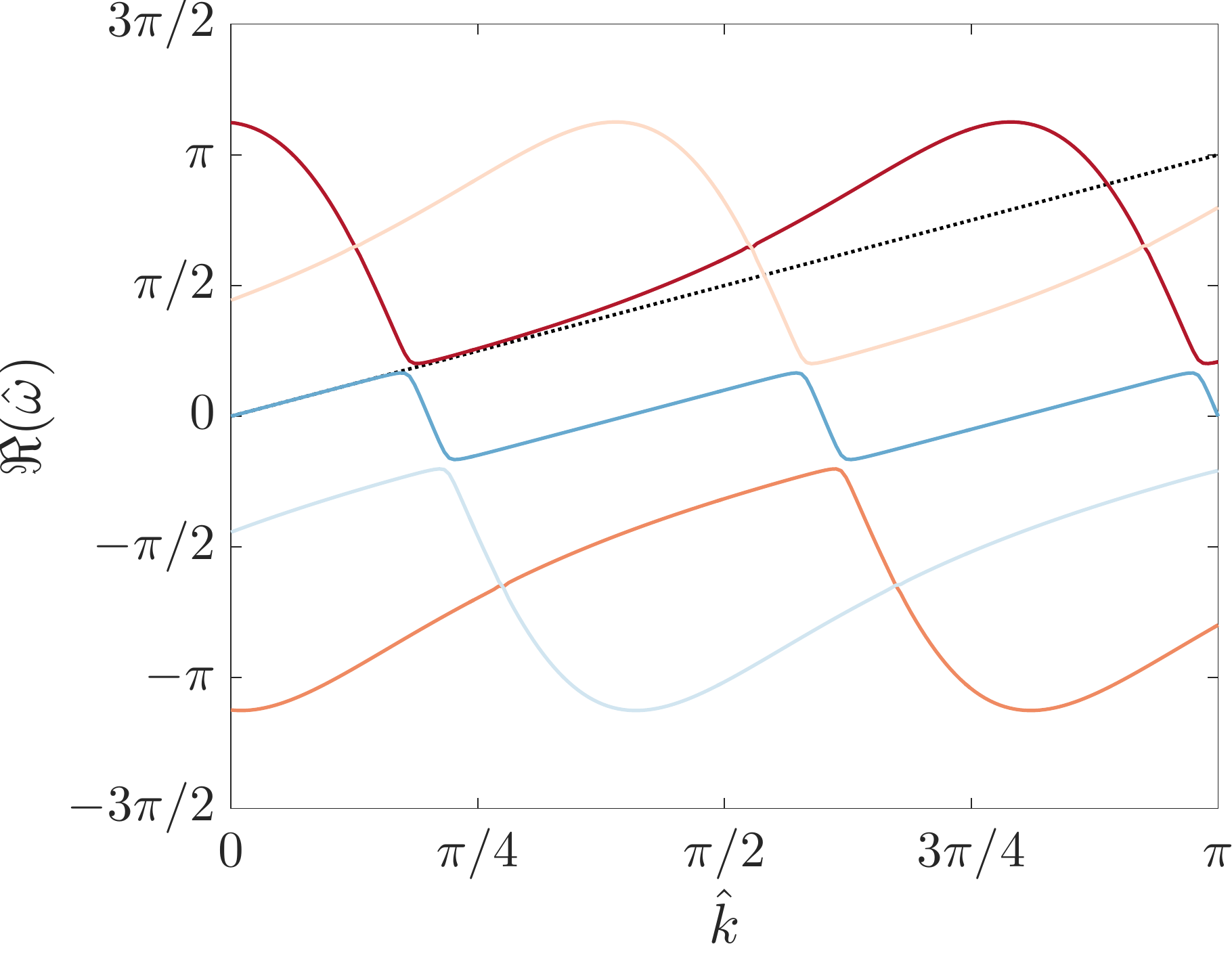}
			\caption{$\mathbf{H}_4 = [0.77,-0.52]$}
			\label{fig:FRGLS4_f1_diss_c}
		\end{subfigure}
		\caption{Linear advection dispersion relations for central interface flux FR, $p=4$, with the dispersion shown for all modes.}
		\label{fig:FRGLS_disp_diss}
	\end{figure}
	
	Lastly for the advection study, the semi-discrete dispersion relation for DG and $\mathbf{H}_4 =[0.77,-0.52]$ with central interfaces are included for completeness. It can be seen that the regions in DG where the dispersion relation is discontinuous has been greatly reduced, from \cite{Watkins2016,Asthana2017}, it can be understood that the smaller regions of discontinuity in the dispersion relation, leads to less energy transfer to spurious modes. 
	
   To study the general trend in the CFL limit of the combined advection diffusion equation we use the case of $a=10$ and $\nu=1$, which was similarly used in the investigation by Watkins~\etal~\cite{Watkins2016}.
   \begin{figure}
		\centering
		\begin{subfigure}[b]{0.45\linewidth}
			\centering
			\includegraphics[width=\linewidth]{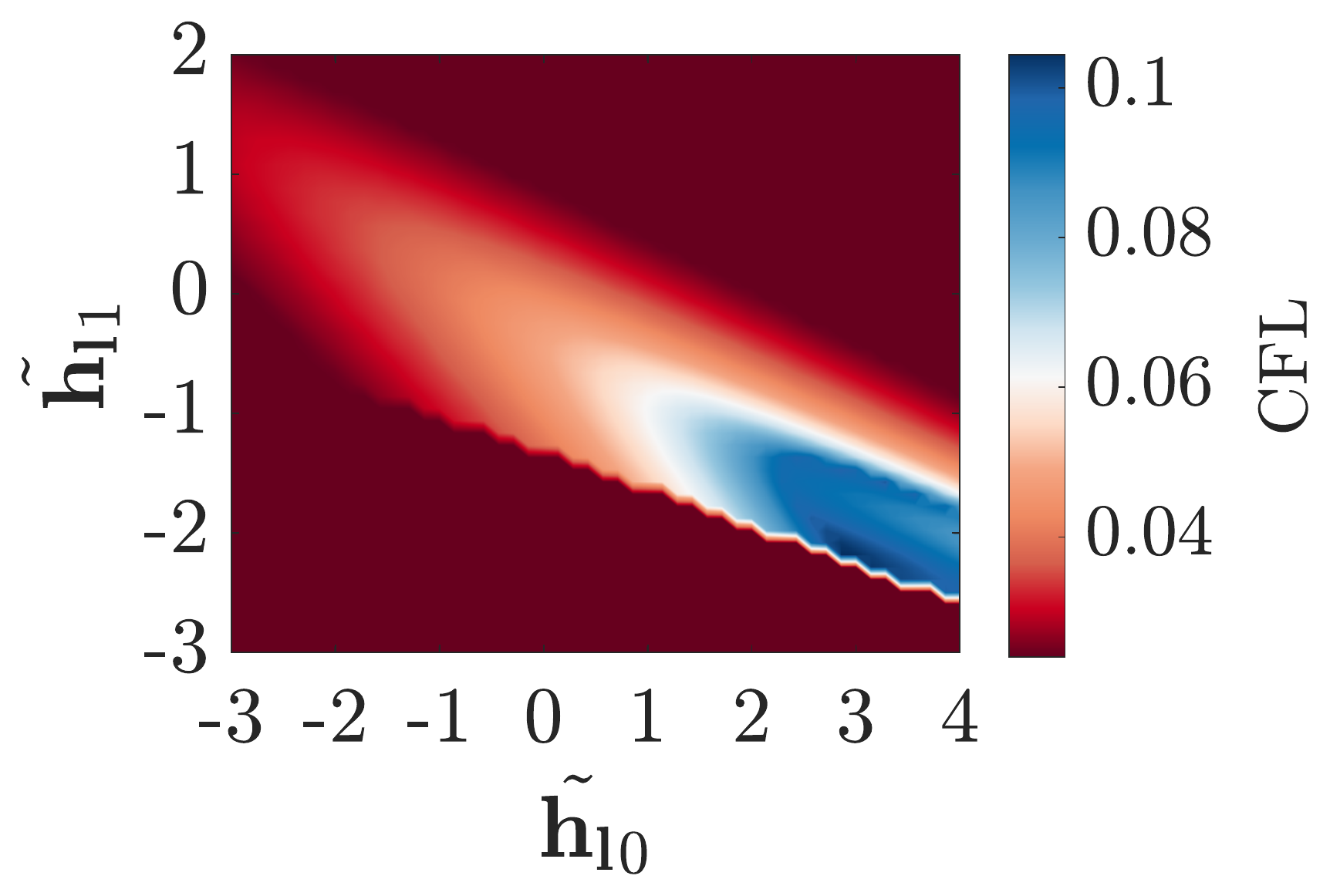}
			\caption{$p=4$, $\alpha_a = 0.5$, $\alpha_d = 0.5$}
			\label{fig:FRGLS4_RK44_ad_CFL_central}
		\end{subfigure}
		~
		\begin{subfigure}[b]{0.45\linewidth}
			\centering
			\includegraphics[width=\linewidth]{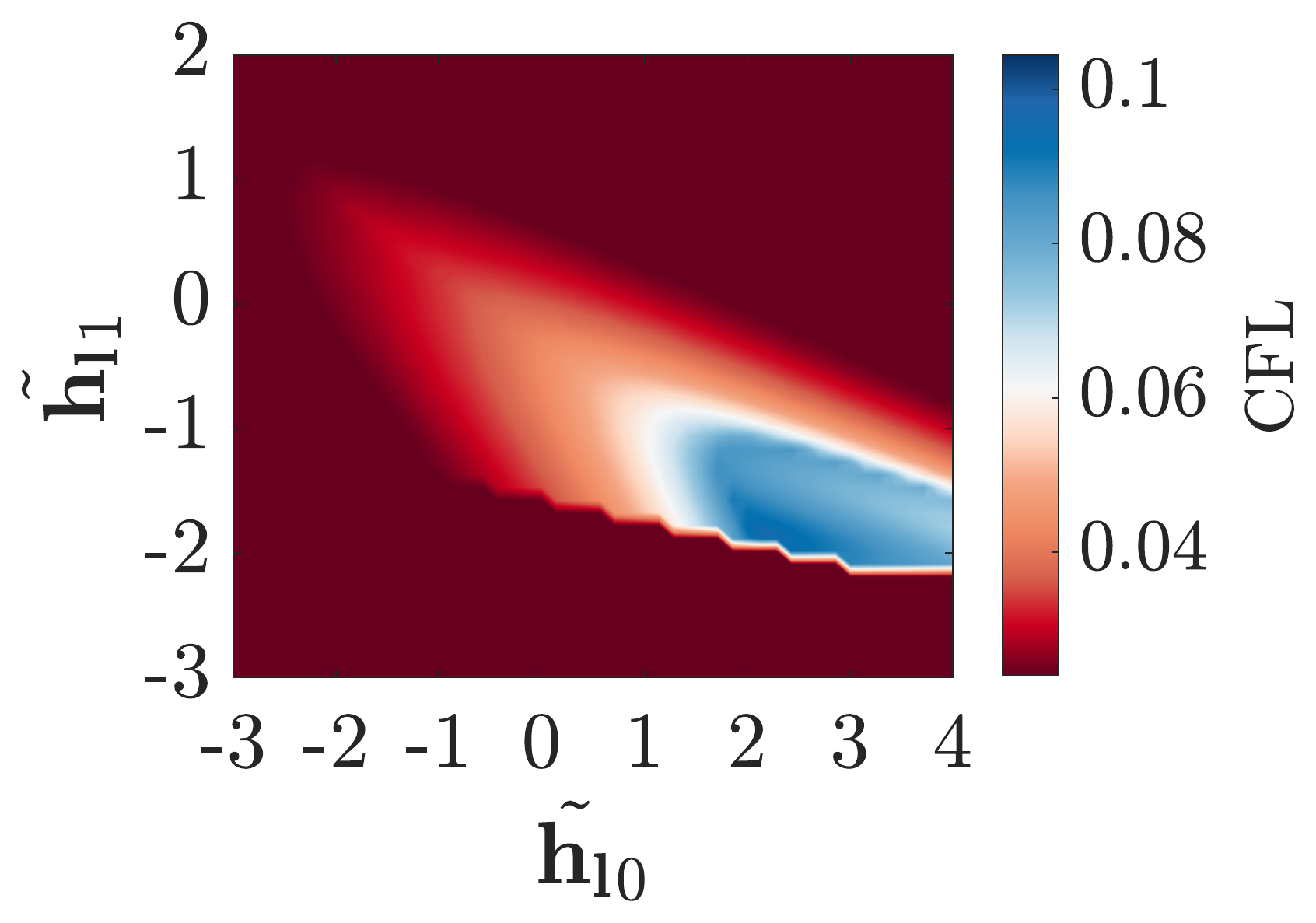}
			\caption{$p=4$, $\alpha_a = 1$, $\alpha_d=0.5$}
			\label{fig:FRGLS4_RK44_ad_CFL_up}
		\end{subfigure}
		\caption{Maximum stable CFL number for GLSFR $[c=10,\nu=1]$, with RK44 temporal integration.}
		\label{fig:GLSFR_ad_CFL}
	\end{figure}
	Figure~\ref{fig:GLSFR_ad_CFL} shows the CFL limits for $p=4$ for various correction functions for the case when the both advection and diffusion interface calculation are central differenced and the case when advection is upwinded and diffusion is central differenced. In both cases it can be seen that the stability profile is of a different character than pure advection, with higher values of $\hat{\mathbf{h_l}}_0$ being stable due to the additional stability brought by the diffusion terms. The motivation of this test is the understanding it can bring to later tests. With the main observation being the slant to the CFL profile, although this will be affected by the relative scale of the advection-diffusion terms.
	
\section{Taylor-Green Vortex}
	So far the investigation has focused on theoretical numerical analysis, and there is already a plethora of papers which similarly investigate the theoretical behaviour of FR. Therefore, to build on the theoretical work carried out in earlier sections, we will introduce the Talyor-Green vortex test case \cite{Taylor1937,Brachet1983,DeBonis2013}. This is a well known canonical case for the full incompressible Naiver-Stokes equations and is highly important as it exhibits: inviscid, transitional, and turbulent flow regimes. Therefore, if improvements can be made via correction function selection on such a case this will translate well onto an engineering CFD problem. The case is set up as such:
	\begin{align}
		u &= U_0\sin{\bigg(\frac{x}{L}\bigg)}\cos{\bigg(\frac{y}{L}\bigg)}\cos{\bigg(\frac{z}{L}\bigg)} \\
		v &= -U_0\cos{\bigg(\frac{x}{L}\bigg)}\sin{\bigg(\frac{y}{L}\bigg)}\cos{\bigg(\frac{z}{L}\bigg)} \\
		w &= 0 \\
		p &= p_0 + \frac{\rho_0U_0^2}{16}\bigg(\cos{\bigg(\frac{2x}{L}\bigg)} + \cos{\bigg(\frac{2y}{L}\bigg)}\bigg)\bigg(\cos{\bigg(\frac{2z}{L}\bigg)} + 2\bigg) \\
		\rho &= \frac{p}{RT_0}
	\end{align}
	Where we define the case by the non-dimensional parameters as:
	\begin{equation}		
		R_e = \frac{\rho_0U_0L}{\mu}, \quad\quad
		P_r = 0.71 = \frac{\mu\gamma R}{\kappa(\gamma-1)}, \quad\quad
		M_a = 0.08 = \frac{U_0}{\sqrt{\gamma R T_0}} 
	\end{equation}
	with the free variable set as:
	\begin{equation}
		U_0 = 1, \quad\quad \rho_0 = 1, \quad\quad p_0 = 100, \quad\quad R = 1, \quad\quad \gamma = 1.4, \quad\quad L = 1
	\end{equation}
	and the Reynold's and Prandtl's numbers are controlled through setting the dynamic viscosity and thermal conductivity.  

   The aim of this case is to display the impact of changing the correction function by evaluating the approximate time averaged error in the turbulent kinetic energy dissipation by comparison to DNS results. The definition used of turbulent kinetic energy dissipation rate is $-\dx{\epsilon}{t}$, with:
   \begin{equation}
       \epsilon = \frac{1}{2\mathbf{\Omega}}\int_{\mathbf{\Omega}}\rho(u^2 + v^2 + w^2)d\mathbf{\Omega}
   \end{equation}
   
   \begin{figure}[h]
	 	\centering
		\begin{subfigure}[b]{0.45\linewidth}
			\centering
			\includegraphics[width=\linewidth]{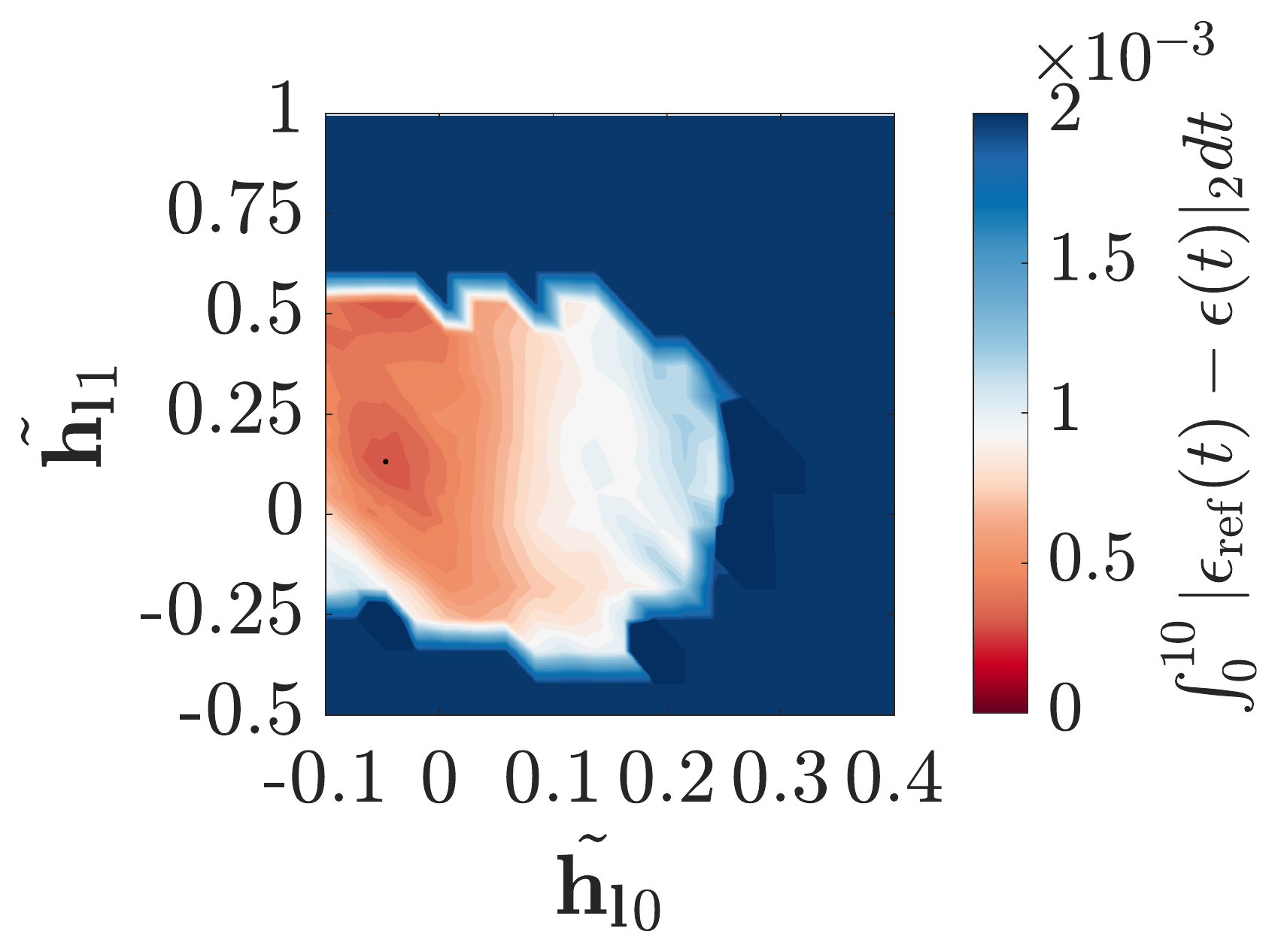}
			\caption{$16^3$ elements}
			\label{fig:FRGLS4_TGV_16_error}
		\end{subfigure}
		~
		\begin{subfigure}[b]{0.45\linewidth}
			\centering
			\includegraphics[width=\linewidth]{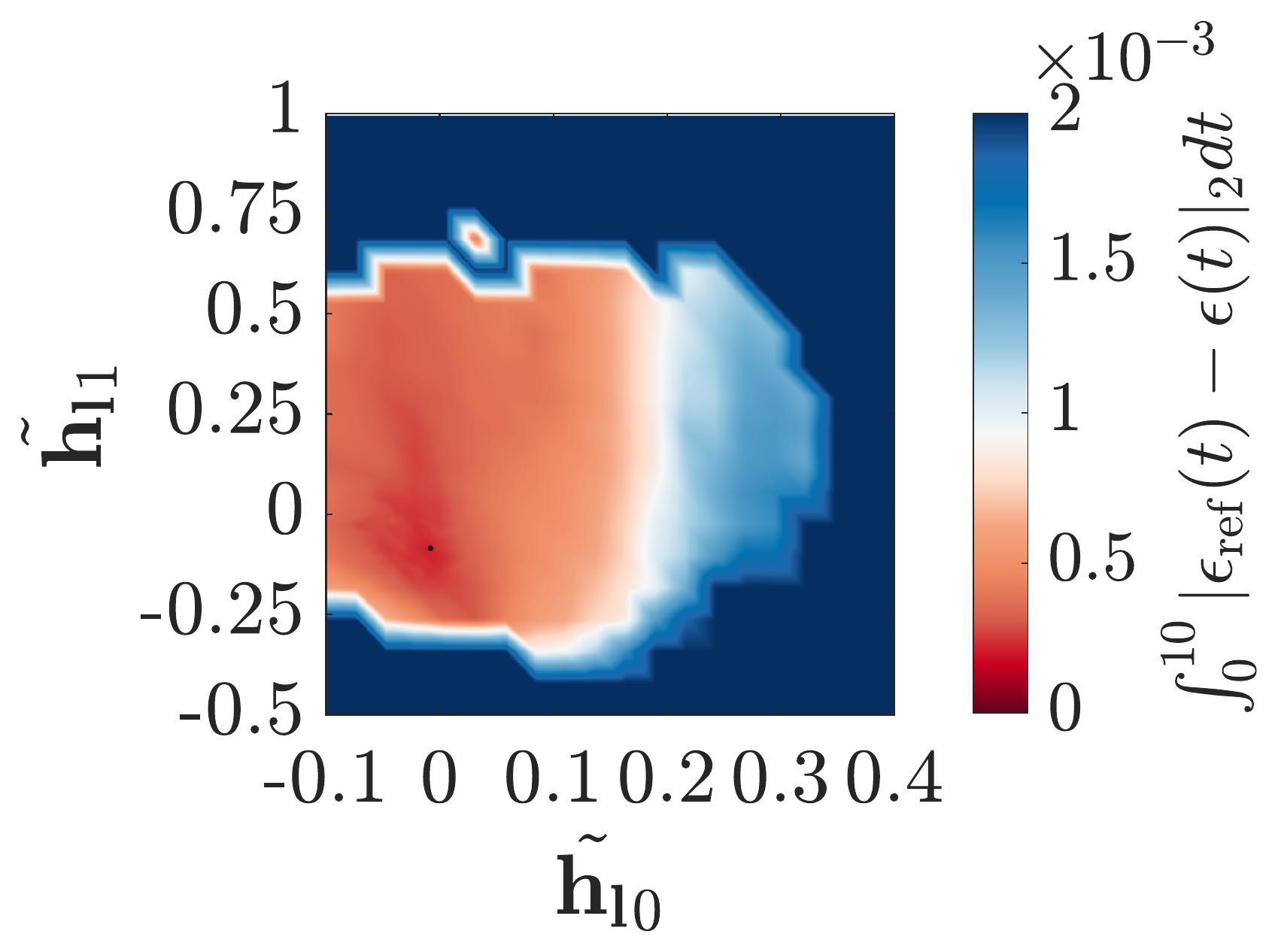}
			\caption{$24^3$ elements}
			\label{fig:FRGLS4_TGV_24_error}
		\end{subfigure}
		\caption{Time integrated error in the kinetic energy dissipation of the Taylor-Green vortex, $R_e=1600$ $P_r = 0.71$ $M_a=0.08$, for FR, $p=4$, with 800 GLSFR correction functions compared to DNS data for each grid resolution. RK44 temporal integration was used with a fixed time step of $1\times10-3$.  A point indicates the correction function tested with the lowest error.}
		\label{fig:FRGLS4_TGV_error}
	\end{figure}
   
   Figure~\ref{fig:FRGLS4_TGV_error} shows the approximate error for two different grid spacing as the correction function is varied with a constant time step of $1\times10^{-3}$. This time step was chosen as it should be sufficient to run with such a time step~\cite{Mastellone2015} and any reduction beyond this would symbolise that a correction function is not sufficiently temporally stable to be considered applicable for industrial problems. As was observed in Fig.~\ref{fig:FRGLS4_RK44_CFL}~\&~\ref{fig:FRGLS4_RK44_CFL_central} as the parameter $\mathbf{h_l}_1$ is increased the temporal stability is reduced and this trend is also seen in the TGV case. There is also large asymmetry between positive and negative values of $\hat{\mathbf{h_l}}_0$ which was also seen in Fig:~\ref{fig:GLSFR_ad_CFL}. However, for Fig.~\ref{fig:GLSFR_ad_CFL} the opposite behaviour is seen, but this was for CFL limit. The implication of a higher CFL limit is a more dissipative scheme, whereas for reduced error a less dissipative scheme is necessary. Hence, the location of the error optimal correction function lying in the left-hand half plane.  A secondary point that may be made about GLSFR correction functions is that for both grid resolution correction functions were found in broadly the same range of $\mathbf{H}_4$ that gave lower error than Nodal DG via FR. Hence, these correction functions are not merely a mathematical exercise but can provide improved performance for real fluid dynamical simulations. 
   
   	\begin{figure}[h]
	 	\centering
		\begin{subfigure}[b]{0.45\linewidth}
			\centering
			\includegraphics[width=\linewidth]{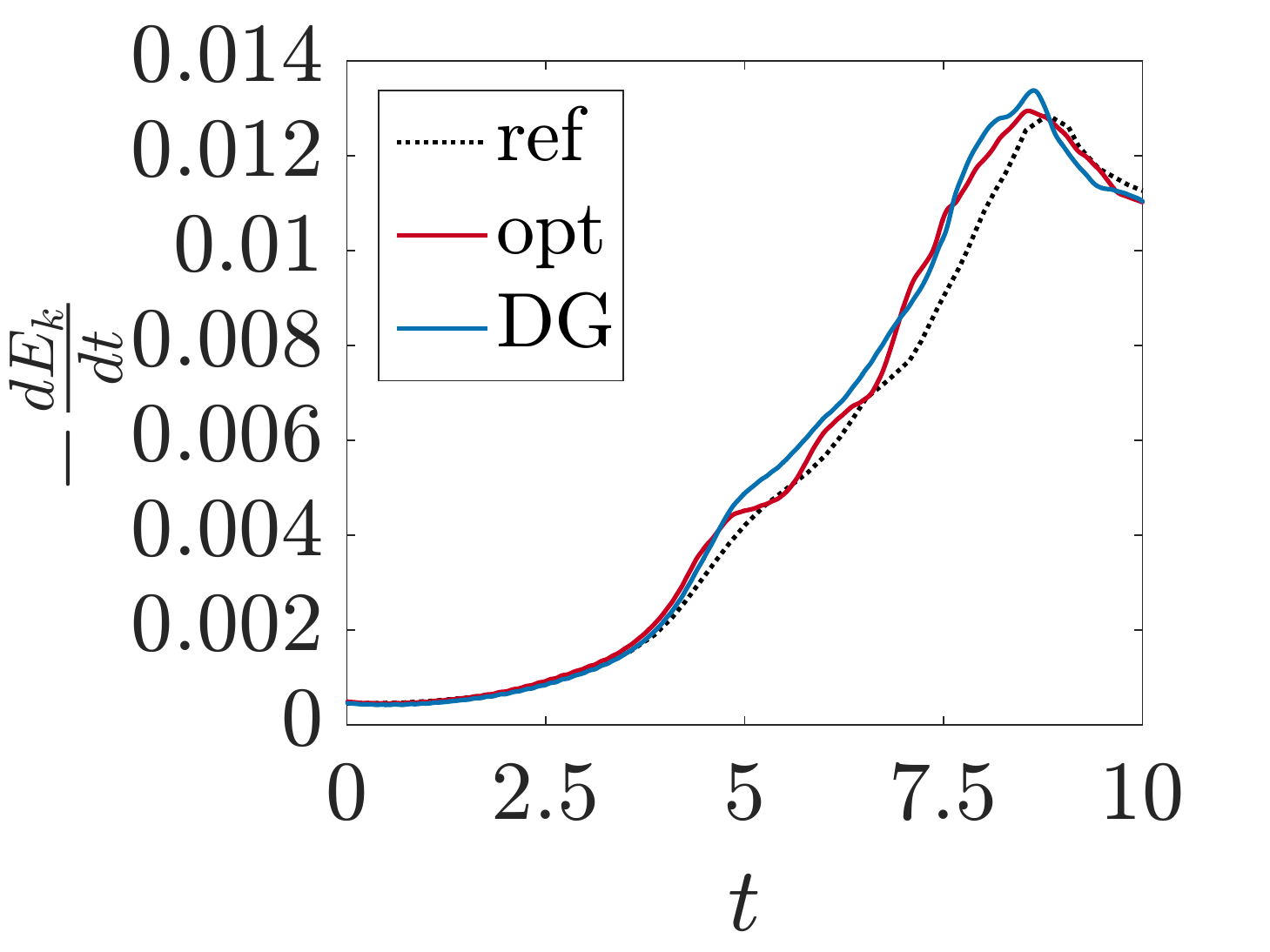}
			\caption{$16^3$ elements}
			\label{fig:FRGLS4_TGV_16_ek}
		\end{subfigure}
		~
		\begin{subfigure}[b]{0.45\linewidth}
			\centering
			\includegraphics[width=\linewidth]{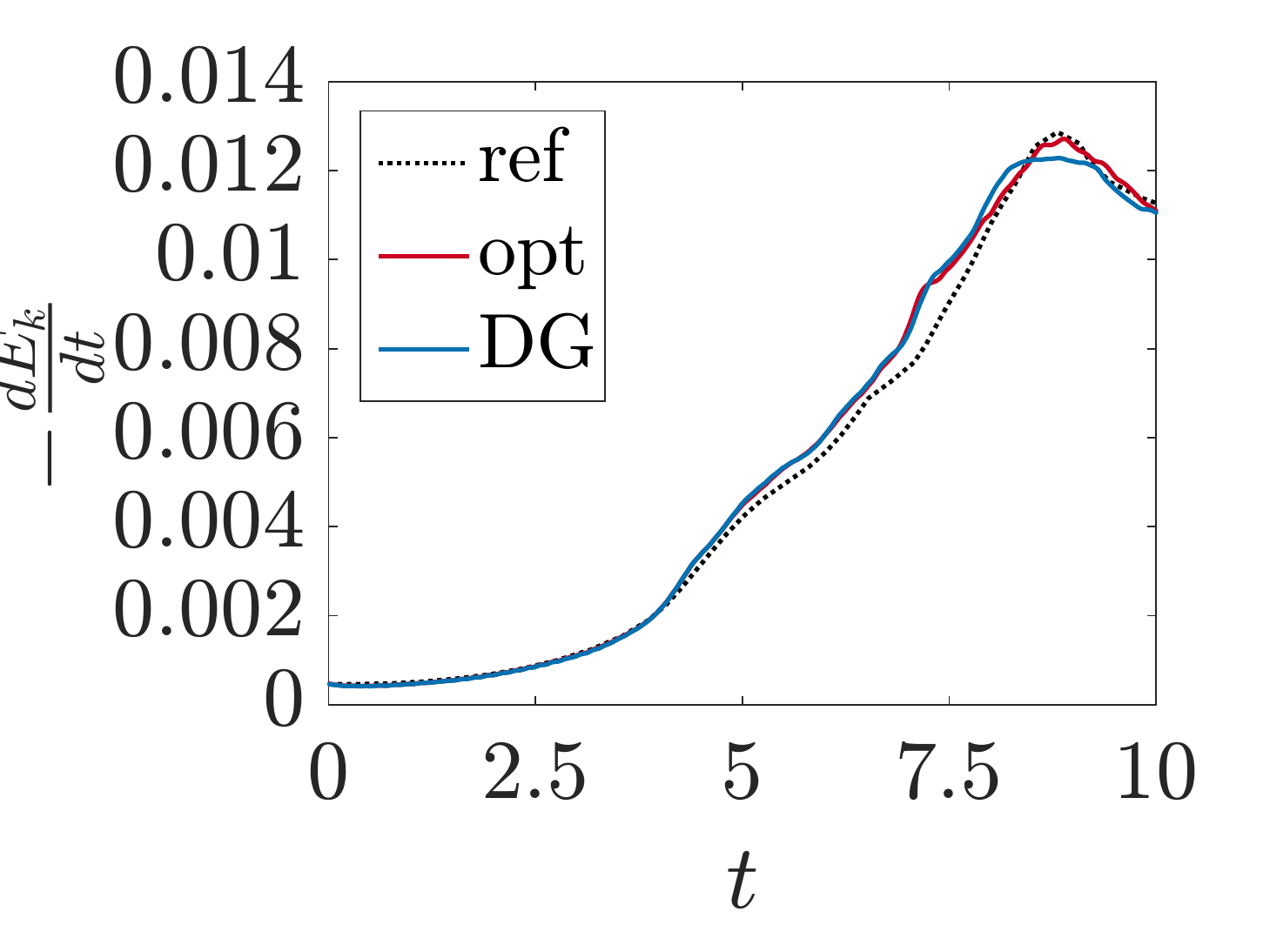}
			\caption{$24^3$ elements}
			\label{fig:FRGLS4_TGV_24_ek}
		\end{subfigure}
		\caption{Optimal GLSFR correction functions, $p=4$, found in Fig.\ref{fig:FRGLS4_TGV_error} compared to Nodal DG via FR, with reference data taken from DeBonis~\cite{DeBonis2013}. RK44 temporal integration was used with a fixed time step of $1\times10-3$.}
		\label{fig:FRGLS4_TGV_ek}
	\end{figure}
	
	Lastly the kinetic energy dissipation of the optimal and DG correction function is compared for the two grid refinement levels in Fig.\ref{fig:FRGLS4_TGV_ek}. In both case the peak turbulent dissipation is better matched, but in both cases there is a region at $t\approx7.5$ in which correction function tuning appear to make only modest improvement. This is the time at which the rate of change of dissipation is highest and the flow can be said to be transitioning to turbulence with turbulent structures forming. The factor can be seen to have the greatest impact on improving the simulation of the physics at this time is increased gird resolution.
   
\section{Conclusions}
\begin{figure}
		\centering
		\includegraphics[width=0.45\linewidth]{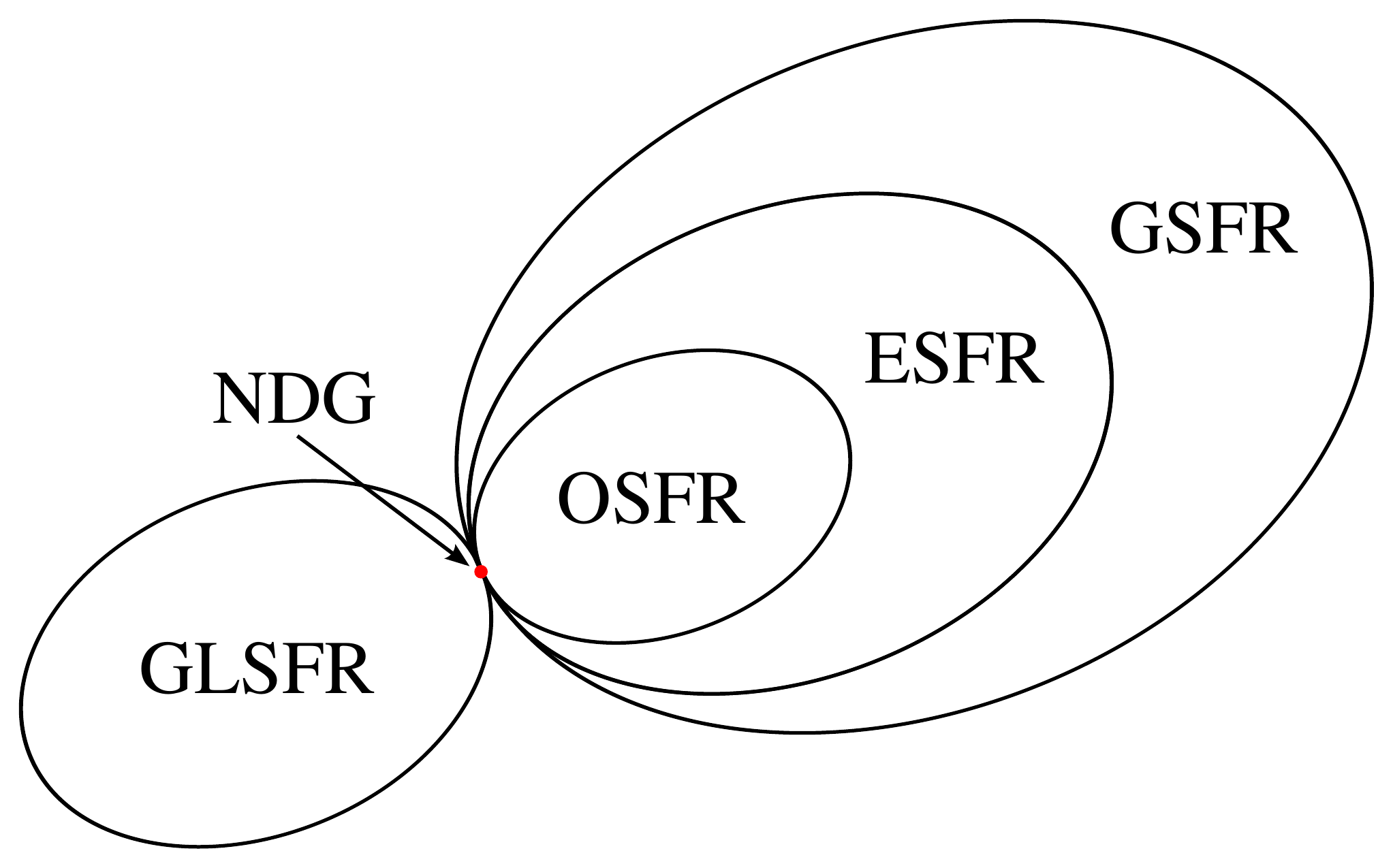}
		\caption{Diagram showing intersection of sets of correction functions and where Nodal DG is located.}
		\label{fig:FR_corr_vd}
	\end{figure}

	It has been shown that through using the Lebesgue norm, a new energy stable family of correction functions can be derived that can be very arbitrary. Through the derivation, it is apparent that this familiy of correction functions extends the current set of stable correction functions, with the only union between the sets of correction being at Nodal DG, Fig.\ref{fig:FR_corr_vd}.

	Von Neumann analysis for both the linear advection and linear advection-diffusion was then presented for this new set of corrections and a region of temporally stable correction function was found to exist for both upwinded and centrally difference interfaces, when Runge-Kutta temporal integration is used. Through exploration of the dispersion and dissipation characteristics of some of the correction function, it was shown that GLSFR may be able to provide improved performance for implicit LES calculations. Lastly, a series of Taylor-Green vortex test are performed for two grid spacing and it was found that in both case a similar region of corrections gave good performance for minimising the error in the turbulent kinetic energy dissipation. In both cases correction functions were found that reduced the error when compared to DG. Therefore, it is concluded that GLSFR may be able improve the accuracy of FR when applied to real fluid-dynamical problems.
	
\section*{Acknowledgements}
The support of the Engineering and Physical Sciences Research Council of the United Kingdom is gratefully acknowledged. The author would like to acknowledge the useful discussions had with Rob Watson.

\bibliographystyle{aiaa}
\bibliography{library}

\begin{thebibliography}{10}
\newcommand{\enquote}[1]{``#1''}

\bibitem{Harten1987}
Harten, A., Engquist, B., Osher, S., and Chakravarthy, S.~R.,
  \enquote{{Uniformly High-Order Accurate Essentially Non-Oscillatory Schemes,
  III},} {\em Journal of Computational Physics\/}, Vol.~71, No.~2, 1987,
  pp.~231--303.

\bibitem{Liu1994}
Liu, X.-D., Osher, S., and Chan, T., \enquote{{Weighted Eessentially
  Non-Oscillatory Schemes},} {\em Journal of Computational Physics\/},
  Vol.~115, No.~1, 1994, pp.~200--212.

\bibitem{Reed1973}
Reed, W.~H. and Hill, T.~R., \enquote{{Triangular Mesh Methods for the Neutron
  Transport Equation},} {\em Los Alamos Report LA-UR-73-479\/}, , No. 836,
  1973, pp.~10.

\bibitem{Cockburn2000}
Cockburn, B., Karniadakis, G.~E., and Shu, C.-W., \enquote{{The Development of
  Discontinuous Galerkin Methods},} {\em Discontinuous Galerkin Methods
  (Newport, RI, 1999)\/}, Vol.~11, chap.~1, Springer-Verlag, Berlin, 2000, pp.
  3--50.

\bibitem{Karniadakis2013}
Karniadakis, G.~E. and Sherwin, S., {\em {Spectral/hp Element Methods for
  Computational Fluid Dynamics}\/}, Oxford University Press, Oxford, 1st ed.,
  2013.

\bibitem{Wang2002}
Wang, Z., \enquote{{Spectral (Finite) Volume Method for Conservation Laws on
  Unstructured Grids. Basic Formulation: Basic Formulation},} {\em Journal of
  Computational Physics\/}, Vol.~178, No.~1, 2002, pp.~210--251.

\bibitem{Gao2013a}
Gao, H. and Wang, Z.~J., \enquote{{A Conservative Correction Procedure via
  Reconstruction Formulation with the Chain-Rule Divergence Evaluation},} {\em
  Journal of Computational Physics\/}, Vol.~232, No.~1, 2013, pp.~7--13.

\bibitem{Wang2009}
Wang, Z.~J. and Gao, H., \enquote{{A unifying lifting collocation penalty
  formulation including the discontinuous Galerkin, spectral volume/difference
  methods for conservation laws on mixed grids},} {\em Journal of Computational
  Physics\/}, Vol.~228, No.~21, 2009, pp.~8161--8186.

\bibitem{Huynh2007}
Huynh, H.~T., \enquote{{A Flux Reconstruction Approach to High-Order Schemes
  Including Discontinuous Galerkin Methods},} {\em 18th AIAA Computational
  Fluid Dynamics Conference\/}, Vol. 2007-4079, 2007, pp. 1--42.

\bibitem{Huynh2009}
Huynh, H.~T., \enquote{{A Flux Reconstruction Approach to High-Order Schemes
  Including Discontinuous Galerkin for Diffusion},} {\em 47th AIAA Aerospace
  Science Meeting\/}, No. January in Fluid Dynamics and Co-located Conferences,
  American Institute of Aeronautics and Astronautics, jun 2009, pp. 1--34.

\bibitem{Williams2011}
Williams, D.~M., Castonguay, P., Vincent, P.~E., and Jameson, A., \enquote{{An
  Extension of Energy Stable Flux Reconstruction to Unsteady , Non-linear ,
  Viscous Problems on Mixed Grids},} {\em Fluid Dynamics\/}, , No. June, 2011,
  pp.~1--37.

\bibitem{Huynh2011}
Huynh, H.~T., \enquote{{High-Order Methods Including Discontinuous Galerkin by
  Reconstructions on Triangular Meshes},} {\em 49th AIAA Aerospace Sciences
  Meeting\/}, No. January, 2011, pp. 1--28.

\bibitem{Castonguay2012a}
Castonguay, P., Vincent, P.~E., and Jameson, A., \enquote{{A new class of
  high-order energy stable flux reconstruction schemes for triangular
  elements},} {\em Journal of Scientific Computing\/}, Vol.~51, No.~1, 2012,
  pp.~224--256.

\bibitem{Williams2014}
Williams, D.~M. and Jameson, A., \enquote{{Energy stable flux reconstruction
  schemes for advection-diffusion problems on tetrahedra},} {\em Journal of
  Scientific Computing\/}, Vol.~59, No.~3, 2014, pp.~721--759.

\bibitem{Sheshadri2016a}
Sheshadri, A. and Jameson, A., \enquote{{On the Stability of the Flux
  Reconstruction Schemes on Quadrilateral Elements for the Linear Advection
  Equation},} {\em Journal of Scientific Computing\/}, Vol.~67, No.~2, 2016,
  pp.~769--790.

\bibitem{Vincent2011}
Vincent, P.~E., Castonguay, P., and Jameson, A., \enquote{{Insights From von
  Neumann Analysis of High-Order Flux Reconstruction Schemes},} {\em Journal of
  Computational Physics\/}, Vol.~230, No.~22, 2011, pp.~8134--8154.

\bibitem{Vermeire2016}
Vermeire, B. and Vincent, P., \enquote{{On the properties of energy stable flux
  reconstruction schemes for implicit large eddy simulation},} {\em Journal of
  Computational Physics\/}, , No. September, 2016.

\bibitem{Vincent2010}
Vincent, P.~E., Castonguay, P., and Jameson, A., \enquote{{A New Class of
  High-Order Energy Stable Flux Reconstruction Schemes},} {\em Journal of
  Scientific Computing\/}, Vol.~47, No.~1, 2010, pp.~50--72.

\bibitem{Vincent2015}
Vincent, P.~E., Farrington, A.~M., Witherden, F.~D., and Jameson, A.,
  \enquote{{An Extended Range of Stable-Symmetric-Conservative Flux
  Reconstruction Correction Functions},} {\em Computer Methods in Applied
  Mechanics and Engineering\/}, Vol.~296, 2015, pp.~248--272.

\bibitem{Trojak2018}
Trojak, W., \enquote{{Generalised Sobolev Stable Flux Reconstruction},}
  Vol.~Arxiv:1804, No. 04714, 2018, pp.~1--19.

\bibitem{Jameson2012}
Jameson, A., Vincent, P.~E., and Castonguay, P., \enquote{{On The Non-Linear
  Stability of Flux Reconstruction Schemes},} {\em Journal of Scientific
  Computing\/}, Vol.~50, No.~2, 2012, pp.~434--445.

\bibitem{Toro2009}
Toro, E., {\em {Riemann Solvers and Numerical Methods for Fluid Dynamics - A
  Practical Introduction}\/}, Springer-Verlag Berlin Heidelberg, Dordrecht
  Berlin Heidelberg London New York, 3rd ed., 2009.

\bibitem{Kennedy2000}
Kennedy, C.~A., Carpenter, M.~H., and Lewis, R.~M., \enquote{{Low-Storage,
  Explicit Runge-Kutta Schemes for the Compressible Navier-Stokes Equations},}
  {\em Applied Numerical Mathematics\/}, Vol.~35, No.~1, 2000, pp.~177--219.

\bibitem{Holdeman1970}
Holdeman, J.~T., \enquote{{Legendre Polynomial Expansions of Hypergeometric
  Functions with Applications},} {\em Journal of Mathematical Physics\/},
  Vol.~11, No.~1, 1970, pp.~114--117.

\bibitem{Trojak2017}
Trojak, W., Watson, R., and Tucker, P.~G., \enquote{{High Order Flux
  Reconstruction on Stretched and Warped Meshes},} {\em AIAA SciTech, 55th AIAA
  Aerospace Sciences Meeting\/}, Grapevine Texas, 2017, pp. 1--12.

\bibitem{Lele1992}
Lele, S.~K., \enquote{{Compact Finite Difference Schemes With Spectral-Like
  Resolution},} {\em Journal of Computational Physics\/}, Vol.~103, No.~1,
  1992, pp.~16--42.

\bibitem{Hu1999}
Hu, F.~Q., Hussaini, M., and Rasetarinera, P., \enquote{{An Analysis of the
  Discontinuous Galerkin Method for Wave Propagation Problems},} {\em Journal
  of Computational Physics\/}, Vol.~151, No.~2, 1999, pp.~921--946.

\bibitem{Hesthaven2008}
Hesthaven, J.~S. and Kirby, R.~M., \enquote{{Filtering in Legendre spectral
  methods},} {\em Mathematics of Computation\/}, Vol.~77, No. 263, 2008,
  pp.~1425--1452.

\bibitem{Williams2014a}
Castonguay, P., Williams, D.~M., Vincent, P.~E., and Jameson, A.,
  \enquote{{Energy Stable Flux Reconstruction Schemes for Advection-Diffusion
  Problems},} {\em Computer Methods in Applied Mechanics and Engineering\/},
  Vol.~267, No.~1, 2013, pp.~400--417.

\bibitem{Castonguay2012}
Castonguay, P., {\em {High-order Energy Stable Flux Reconstruction Schemes For
  Fluid Flow Simulations on Unstructured Grids}\/}, Ph.D. thesis, Stanford
  University, 2012.

\bibitem{Watkins2016}
Watkins, J., Asthana, K., and Jameson, A., \enquote{{A Numerical Analysis of
  the Nodal Discontinuous Galerkin Scheme Via Flux Reconstruction for the
  Advection-Diffusion Equation},} {\em Computers and Fluids\/}, Vol.~139,
  No.~1, 2016, pp.~233--247.

\bibitem{Trefethen1994}
Trefethen, L.~N., {\em {Finite Difference and Spectral Methods for Ordinary and
  Partial Differential Equations}\/}, Ithaca, New York, 1st ed., 1994.

\bibitem{Isaacson1994}
Isaacson, E. and Keller, H.~B., {\em {Analysis of Numerical Methods}\/}, John
  Wilery {\&} Sons Ltd., New York, 2nd ed., 1994.

\bibitem{Vincent2017}
Vincent, P.~E., Witherden, F., Vermeire, B., Park, J.~S., and Iyer, A.,
  \enquote{{Towards Green Aviation with Python at Petascale},} {\em
  International Conference for High Performance Computing, Networking, Storage
  and Analysis, SC\/}, No. November, 2017, pp. 1--11.

\bibitem{Asthana2017}
Asthana, K., Watkins, J., and Jameson, A., \enquote{{On Consistency and Rate of
  Convergence of Flux Reconstruction For Time-Dependent Problems},} {\em
  Journal of Computational Physics\/}, Vol.~334, 2017, pp.~367--391.

\bibitem{Taylor1937}
Taylor, G.~I. and Green, A.~E., \enquote{{Mechanism of the Production of Small
  Eddies from Large Ones},} {\em Proceedings of the Royal Society A:
  Mathematical, Physical and Engineering Sciences\/}, Vol.~158, No. 895, 1937,
  pp.~499--521.

\bibitem{Brachet1983}
Brachert, M.~E., Orszag, S.~A., Nickel, B.~G., Morf, R.~H., and Frisch, U.,
  \enquote{{Small-Scale Structure of the Taylor-Green Vortex},} {\em Journal of
  Fluid Mechanics\/}, Vol.~130, 1983, pp.~411--452.

\bibitem{DeBonis2013}
DeBonis, J., \enquote{{Solutions of the Taylor-Green Vortex Problem Using
  High-Resolution Explicit Finite Difference Methods},} {\em 51st AIAA
  Aerospace Sciences Meeting including the New Horizons Forum and Aerospace
  Exposition\/}, , No. February 2013, 2013.

\bibitem{Mastellone2015}
Mastellone, A., Capuano, F., Benedetto, S., and Cutrone, L., \enquote{{Problem
  C3 . 5 Direct Numerical Simulation of the Taylor-Green Vortex at Re = 1600},}
  Tech. rep., 2015.

\end{thebibliography}

\newpage
\appendix

\section{GLSFR Generating Pseudo-Code}
   We will now present a short piece of pseudo-code to calculate the left correction functions Legendre polynomial coefficients, $\tilde{\mathbf{h_l}}$. In this example we use $\%$ to mean the modulus or remainder operator, \emph{i.e.} $a\% b = a \mod b$. Here the input is an array $\mathbf{q}$ with $p-2$ entries that control the shape of $h_L$.

	\begin{algorithm}
		\caption{Calculate left GLSFR correction Legendre coefficients}
		\begin{algorithmic}
			\REQUIRE $p,q[p-2]$
			\STATE $\tilde{h_L}[p-1] \leftarrow 0$
			\STATE $\tilde{h_L}[p-2] \leftarrow 0$
			\FOR{$0\leqslant i \leqslant p-3$} 
				\STATE $\tilde{h_L}[i] \leftarrow q[i]$
				\STATE $\tilde{h_L}[p-1] \leftarrow \tilde{h_L}[p-1] - (i+p\%2)q[i]$
				\STATE $\tilde{h_L}[p-2] \leftarrow \tilde{h_L}[p-2] - (i+p-1\%2)q[i]$
			\ENDFOR
			\STATE $\tilde{h_L}[p] \leftarrow 0.5(-1)^p$
			\STATE $\tilde{h_L}[p+1] \leftarrow 0.5(-1)^{p+1}$
			\RETURN $\tilde{h_L}$
		\end{algorithmic}
	\end{algorithm}

\section{Nomenclature}
		\begin{tabbing}
		  XXXXXX \= \kill% this line sets tab stop  
		  	\textit{Roman}\\
		  	$a_p$ \> $(2p)!/(2^p(p!)^2)$ \\
			$c(k)$ \> modified phase velocity at wavenumber $k$ \\
			$\mathbf{B}_{+2}$ \> diffusion second downwind cell FR matrix \\
			$\mathbf{B}_{+1}$ \> diffusion first downwind cell FR matrix \\
			$\mathbf{B}_0$ \> diffusion centre cell FR matrix \\
			$\mathbf{B}_{-1}$ \> diffusion first upwind cell FR matrix \\
			$\mathbf{B}_{-2}$ \> diffusion second upwind cell FR matrix \\
			$\mathbf{C}_{+1}$ \> advection downwind cell FR matrix \\
			$\mathbf{C}_0$ \> advection centre cell FR matrix \\
			$\mathbf{C}_{-1}$ \> advection upwind cell FR matrix \\
			$\mathbf{D}$ \> first derivative matrix \\
			$f$ \> flux variable in physical domain \\
			$g_L \:\mathrm{\&}\: g_R$ \> left and right correction function gradients\\
			$h_L \:\mathrm{\&}\: h_R$ \> left and right correction functions\\
			$J_i$ \> $i^{\mathrm{th}}$ cell Jacobian\\
			$k$ \> wavenumber \\
			$k_{nq}$ \> solution point Nyquist wavenumber, $(p+1)/\delta_j$\\
			$\hat{k}$ \> $k_{nq}$ normalised wavenumber, $[0,\pi]$ \\
			$\mathbf{K}$ \> ESFR correction matrix \\ 
			$l_i$ \> $i^{\mathrm{th}}$ Lagrange basis function \\
			$\mathbf{L}_p$ \> $p^{\mathrm{th}}$ order GSFR correction matrix \\
			$\mathbf{M}$ \> polynomial basis mass matrix \\ 
			$M_a$ \> Mach number \\
			$p$ \> solution polynomial order \\
			$P_r$ \> Prandlt number \\
			$\mathbf{Q}_a$ \> FR spatial advection discretisation operator matrix \\
			$\mathbf{Q}_d$ \> FR spatial difussion discretisation operator matrix \\
			$\mathbf{Q}_{ad}$ \> FR spatial advection-diffusion discretisation operator matrix \\
			$\mathbf{R}$ \> FR spatial-temporal update matrix \\
			$R_e$ \> Reynolds number \\
			$u$ \> conserved variable in the physical domain \\
			
			\\ \textit{Greek}\\
			$\alpha$ \> interface upwinding ratio ($\alpha = 1 \Rightarrow$ upwinded, $\alpha = 0.5 \Rightarrow$ central) \\
			$\alpha_a$ \> advection interface upwinding ratio \\
			$\alpha_d$ \> diffusion interface upwinding ratio  \\
			$\delta_j$ \> mesh spacing, $x_j-x_{j-1}$\\
			$\epsilon$ \> domain average kinetic energy \\
			$\iota$ \> OSFR correction function parameter \\
			$\iota_i$ \> $i^{\mathrm{th}}$ GSFR correction function parameter \\
			$\nu$ \> diffusion constant \\			
			$\xi$ \> transformed spatial variable \\
			$\rho(\mathbf{A})$ \> spectral radius of $\mathbf{A}$ \\ 
			$\tau$ \> time step \\
			$\psi_i$ \> $i^{\mathrm{th}}$ Legendre polynomial of the first kind \\ 
			$\mathbf{\Omega}$ \> solution domain \\
			$\mathbf{\Omega}_n$ \> $n^{\mathrm{th}}$ solution sub-domain \\
			$\hat{\mathbf{\Omega}}$ \> reference sub-domain\\
			
			\\ \textit{Subscript}\\
			$\mathrm{\bullet}_L$ \> variable at left of cell\\
			$\mathrm{\bullet}_R$ \> variable at right of cell\\
			
			\\ \textit{Superscript}\\
			$\mathrm{\bullet}^I$ \> common value at interface\\
			$\mathrm{\bullet}^T$ \> vector or matrix transpose \\
			$\mathrm{\bullet}^{\delta}$ \> discontinuous value\\
			$\hat{\mathrm{\bullet}}$ \> variable transformed to reference domain\\
			$\tilde{\mathrm{\bullet}}$ \> variable transformed to Legendre basis \\
			
			\\ \textit{Abbreviations}\\
			OSFR \> Original Stable Flux Reconstruction \cite{Vincent2010}\\
			ESFR \> Extended Stable Flux Reconstruction \cite{Vincent2015}\\
			GSFR \> Generalised Sobolev stable Flux Reconstruction \cite{Trojak2018}\\
			GLSFR \> Generalised Lebesgue Stable Flux Reconstruction\\
			
		\end{tabbing}

\end{document}